\documentclass[10pt]{article}
\usepackage[pdftex]{graphicx}
\usepackage{amsmath, amssymb}
\usepackage{amsthm} 
\usepackage[all]{xy}
\usepackage{tikz}

\def\Sch{\text{\rm Sch}}
\def\F{\cal F}
\def\L{{\cal L}}

\def\bgn{\begin}

\def\CL{\text{\rm CL}}

\def\GL{\text{\rm GL}}
\def\E{\bold E}
\def\J{{\cal J}}
\def\L{{\cal L}}

\def\1{{[1]}}
\def\2{{[2]}}
\def\3{{[3]}}
\def\({\left(}
\def\){\right)}
\def\s-circ{\,{\scriptstyle{\circ}}\,}
\def\<<{<\negthinspace \negthinspace<}
\def\Ad{\text{\rm Ad}}
\def\ad{\text{\rm ad}}

\def\even{\text{\rm even}}
\def\bgn{\begin}

\def\endaln{\end{align}}
\def\cou{\ss\text{\rm co}}

\def\<{<\negthinspace \negthinspace <}

\def\t{\theta}

\def\({\left(}
\def\){\right)}
\def\Im{\text{\rm Im}}
\def\Re{\ss\text{\rm Re}}
\def\[{\big[\neg\big[}
\def\]{\big]\neg\big]}
\def\al{\al}

\def\M{{\cal M}}

\def\tr{\text{\rm tr}}
\def\a{\alpha}
\def\b{\beta}

\def\e{\varepsilon}

\def\Gam{\Gamma}

\def\del{\delta}
\def\lam{\lambda}
\def\ome{\omega}
\def\Ome{\Omega}
\def\sig{\sigma}

\def\E{{\cal E}}

\def\Diff{\text{\rm Diff}_0}

\def\R{\Bbb R}
\def\C{\Bbb C}

\def\Z{\Bbb Z}

\def\M{\frak M}

\def\X{{\cal X}}
\def\w{\wedge}
\def\({\left(}
\def\){\right)}

\def\neg{\negthinspace}
\def\h{\hat}

\def\wtil{\widetilde}

\def\ol{\overline}
\def\pa{\partial}

\def\ran{\rangle} 
\def\lan{\langle}

\def\ss{\scriptscriptstyle}
\def\trian{\triangle}

\def\Ker{\text{\rm Ker}}
\def\arrow{\longrightarrow}

\def\:{\, :\,}
\def\CL{\text{\rm CL}}
\def\TT{T\oplus T^*}
\def\complex{generalized complex }
\def\K\"ahler{generalized K\"ahler}
\def\vol{\text{\rm vol}}

\def\10{\displaystyle L^{10}}
\def\2{\displaystyle L^2}
\def\c0{\displaystyle C^0}

\def\10{\displaystyle L^{10}}
\def\2{\displaystyle L^2}
\def\del{\delta}
\def\del2{\displaystyle L^2_{0,\delta}}
\def\c0{\displaystyle C^0}

\def\del{\delta}

\def\K{{\cal K}}
\def\M-A{\text{\rm Monge-Amp\`ere}}

\def\M-A{\text{\rm Monge-Amp\`ere}}
\def\[{\big[\,}
\def\]{\,\big]}
\def\End{\text{\rm End\,}}
\def\Hom{\text{\rm Hom}}

\def\id{\text{\rm id}}

\def\Ham{\text{\rm Ham}}
\def\B{\mathcal B}
\def\even{\text{\rm even}}
\def\odd{\text{\rm odd}}

\def\L{\text{{\cal L}}}
\def\TT{T_M\oplus T^*_M}

\def\M{{\cal M}}
\def\Diff{\text{\rm Diff}}
\def\ol{\overline}
\def\part{\partial}

\def\tt{\ss T\oplus T^*}

\def\F{{\cal F}}
\def\L{{\cal L}}

\def\ad{\text{\rm ad}}

\def\red{0}

\def\integral{\text{\rm int}}
\def\ham{\frak ham}
\def\Re{\R}

\theoremstyle{plain} 

\theoremstyle{definition}

\begin{document}

\title{Matsushima-Lichnerowicz type theorems of Lie algebras of automorphisms of generalized K\"ahler manifolds\\
of symplectic type
%Insert your title here%\thanks{Grants or other notes
}


\author{Ryushi Goto
      %  Second Author %etc.
}


\maketitle
%\tableofcontents
\begin{abstract}
 In K\"ahler geometry, Fujiki--Donaldson show that the scalar curvature arises as the moment map for Hamiltonian diffeomorphisms. In generalized K\"ahler geometry, one does not have  good notions of Levi-Civita connection and curvature, however there still exists a precise framework for the moment map and the scalar curvature is defined as the moment map \cite{Goto_2016}. Then a fundamental question is to understand the existence or non-existence of generalized K\"ahler structures with constant scalar curvature. In the paper, we study the Lie algebra of automorphisms of a generalized complex manifold $(M,\J).$ 
We assume that $H^1(M,\C)=0.$  
Then we show that the Lie algebra of the automorphisms is a reductive Lie algebra if $(M,\J)$ admits a generalized K\"ahler structure of symplectic type with constant scalar curvature. This is a generalization of Matsushima and Lichnerowicz theorem in K\"ahler geometry. 
We explicitly calculate the Lie algebra of the automorphisms of a generalized complex structure given by a cubic curve on $\C P^2$. Cubic curves are classified into nine cases (see Figure.$1\sim9$ in Section \ref{The Lie alg CP2}).
In the three cases as in Figures.\,\,7, 8 and 9,  the Lie algebra of the automorphisms is not reductive and 
there is an obstruction to the existence of generalized K\"ahler structures of symplectic type with constant scalar curvature
in the three cases. 
We also discuss deformations starting from an ordinary K\"ahler manifold $(X,\ome)$ with constant scalar curvature and show that 
nontrivial generalized K\"ahler structures of symplectic type with constant scalar curvature arise as deformations
if the Lie algebra of automorphisms of $X$ is trivial.
We show the Hessian formula of generalized extremal K\"ahler structures and obtain the decomposition theorem of 
the Lie algebra of the reduced automorphisms of a generalized extremal K\"ahler manifold.

%\keywords{Generalized K\"ahler manifold\and moment map \and Matsushima-Lichnerowicz type theorem}
 %{53D18 \and 53C25\and 53D20}
 %MSC code1 \and MSC code2 \and more}
\end{abstract}
\numberwithin{equation}{section}
\section{Introduction}
%It is a fundamental problem to study automorphism groups of geometric structures. 
It is known that the isometry group of a compact Riemannian manifold is a compact Lie group of finite dimension \cite{MS_1939}
 and the automorphism group Aut$(X)$ of a compact complex manifold $X$ is a finite dimensional complex Lie group \cite{S.Kob_1972}. 
 The reduced automorphisms Aut$_{\red}(X)$ is defined to be the subgroup of Aut$(X)$ which acts trivially on the Albanese torus $H^0(X, \Ome^1)^*/H_1(X, \Z)$.  Matsushima and Lichnerowicz show that 
if a compact complex manifold $X $ admits a K\"ahler metric with constant scalar curvature, 
then the Lie algebra of the reduced automorphisms Aut$_{\red}(X)$ is a reductive Lie algebra which is the complexification of the isometry group \cite{Li_1958}, \cite{Ma_1957}.
A generalized K\"ahler manifold is a successful generalization of the ordinary K\"ahler manifold.
A generalized K\"ahler structure on a manifold is a triple $(g, I,J)$ consisting of a Riemannian metric $g$ compatible 
with two complex structures $I, J$ satisfying the certain conditions, which has an origin in a non-linear sigma model in Mathematical physics.
However, a generalized K\"ahler structure has a natural description using the language of Hitchin's generalized complex geometry \cite{Hi_2003}
, which is a commutative pair $({\cal J}_1, {\cal J}_2)$ of generalized complex structures equipped with a positivity condition \cite{Gua_2003}.
The deformation-stability theorem of generalized K\"ahler structures shows that every holomorphic Poisson structure on a compact K\"ahler manifold gives rise to nontrivial deformations of generalized K\"ahler structures {\cite{Goto_2009}, \cite{Goto_2010}, 
\cite{Goto_2012}, \cite{Goto_2014}, \cite{Goto_2016}, \cite{Gua_2018}, \cite{Hi_2007}}. 
\\
\indent
In our previous paper \cite{Goto_2016}, the scalar curvature of a generalized K\"ahler manifold of symplectic type is defined as the moment map 
for Hamiltonian diffeomorphisms, which is a natural generalization of the moment map framework due to Fujiki and Donaldson in 
K\"ahler Geometry.
In this paper, the existence and non-existence problems of generalized K\"ahler structures with constant scalar curvature are discussed.
In order to obtain results of the non-existence, we define the Lie algebra of automorphisms of a generalized complex manifold 
and we introduce the Lie algebra of the reduced automorphisms of a generalized K\"ahler manifold.
We show that the Lie algebra of reduced automorphisms is a reductive Lie algebra if there exists a generalized K\"ahler structure of symplectic type with constant scalar curvature (see Definition \ref{Generalized Kahler structure of symplectic type} and Section \ref{The Lie algebra of cscGK}). 
This is a generalization of Matsushima-Lichnerowicz theorem, which gives an 
obstruction to the existence.
\\ \indent
 Regarding the existence questions,
 we discuss deformations starting from an ordinary K\"ahler manifold $(X,\ome)$ with constant scalar curvature and show that 
nontrivial generalized K\"ahler structures of symplectic type with constant scalar curvature arise as deformations
if the Lie algebra of automorphisms of $X$ is trivial (see Theorem \ref{existenceresult}) cf. \cite{LS_1994}.

%%%%%%%%%%%%%%%%%%%%%%%%%%%%%%%%%%%%%%%%%%%%%%%%%%%%%%
This paper is organized as follows. 
In Section \ref{GCS GK}, notations and preliminary results on generalized complex structures and generalized K\"ahler structures are explained. 
There are many good references and lecture notes on them (see for instance \cite{Gua_2003}, \cite{Gua_2011}, \cite{Cav_2007}, 
\cite{Hi_2011}, \cite{Goto_2010}).
In Section \ref{GHamDiff}, we introduce generalized Hamiltonian diffeomorphisms for an arbitrary generalized complex structure. In particular, if a generalized complex structure comes from a real symplectic structure and a $b$-filed, then the group of generalized Hamiltonian diffeomorphisms coincides with the group of  Hamiltonian diffeomorphisms twisted by the $b$-filed. In Section \ref{GSCasmoment}, we show that the set of compatible almost generalized complex structures with a fixed $\J_\psi$ is an infinite dimensional K\"ahler manifold. 
If the generalized complex structure $\J_\psi$ is locally given by $d$-closed nondegenerate, pure spinor, then 
there exists a moment map for the action of generalized Hamiltonian diffeomorphisms with respect to $\J_\psi$.
Then we define the scalar curvature as the moment map. 
In Subsection \ref{The Lie algebra and the reduced Lie algebra}, we introduce the Lie algebra $\frak g_\J$ of automorphisms of a generalized complex manifold and the Lie algebra of the reduced automorphisms 
$\frak g_{\red}.$
In Subsection \ref{the rela Lie algebra frak gredR}, we introduce a real Lie algebra $\frak g_\red^{\Re},$ 
which can be regarded as the Lie algebra of generalized isometry group of a generalized K\"ahler manifold.
We show that $\frak g_{\red}^{\Re}$
is always a reductive Lie algebra.
In Subsection \ref{Reductiveity of frak gred}, we show that the condition (\ref{uRe uIm}) implies $\frak g_\red $
is the complexification of $\frak g_{\red}^{\Re},$ which is then reductive.
In Subsection \ref{The structure theorem of the Lie algebra}, we obtain a structure theorem of $\frak g_\J$
(see Theorem \ref{the structure theorem of frakgJ}). In particular, we see that 
$\frak g_\J\cong\frak g_{\red}$ if $H^{1}(M,\C)=0$ (Corollary \ref{cor hodd=0}).
In Subsection \ref{The Lie algebra poisson deformations},
we show if $H^1(M, {\cal O})=0,$ then the Lie algebra $\frak g_{\J_\b}$ is given by 
the Lie algebra of holomorphic vector fields preserving a holomorphic Poisson structure $\b.$
In Section \ref{The Lie algebra of cscGK}, we prove one of our main theorems :
\\ 

{\indent\sc Theorem} \ref{main theorem}.
\,Let $(M, \J)$ be a $2n$ dimensional compact generalized complex manifold. 
We assume that $H^{1}(M,\C)=0.$
If $M$ admits a generalized K\"ahler structure $(\J, \J_\psi)$ of symplectic type with constant 
scalar curvature, then 
the Lie algebra $\frak g_\J$ is a reductive Lie algebra.
\\ \\
\noindent In fact, the condition $H^{1}(M,\C)=0$ implies that $\frak g_\J =\frak g_{\red}$ and the existence of generalized K\"ahler structure of symplectic type with constant scalar curvature 
implies that  $\frak g_{\red}$ is the complexification of $\frak g_{\red}^{\Re}$ which is reductive (Theorem \ref{the Lie algebra of the reduced automorphisms }).
Applying Theorem \ref{main theorem}
 to the most important cases of Poisson deformations, we obtain \\ \\
{\indent\sc Theorem} \ref{main theorem Poisson deformations}.\,
Let $(M, I, \ome)$ be a compact K\"ahler manifold with a holomorphic Poisson structure $\b\neq 0.$
We assume $H^{1}(M, \C)=0.$
We denote by $(M, \J_{\b t}, \J_{\psi_t})$  a generalized K\"ahler manifold given by Poisson deformations.
Then if the scalar curvature $S(\J_{\b t}, \J_{\psi_t})$ is a constant, the Lie algebra of the automorphisms $\frak g_{\J_{\b t}}$ is a reductive Lie algebra. 
\\
\\
In Section \ref{The Lie alg CP2}, 
 on $\C P^2$, a holomorphic Poisson structure $\b$ is given by a section of the anticanonical line bundle 
with the zero locus given by 
a cubic curve.
Cubic curves of $\C P^2$ are classified into nine cases.
An explicit calculation of the Lie algebra $\frak g_{\J_\b}$ is shown for each case.
In Section \ref{DeformcscGK}, 
the results of the existence are discussed by using deformations.
In particular, Del Pezzo surfaces with trivial automorphisms admit generalized K\"ahler structures with constant
scalar curvature.
In Section \ref{GEXKandHessian}, we introduce a generalized extremal K\"ahler manifold and 
calculate the Hessian of the Calabi type functional. We obtain the decomposition of the Lie algebra of automorphisms of a generalized extremal K\"ahler manifold (cf. \cite{AM_2019}, \cite{FuO}, \cite{LiWang_2006}). 

%%%%%%%%%%%%%%%%%%%%%%%%%%%%%%%%%%%%%%%%%%%%%%%

\section{Generalized complex structures and generalized K\"ahler structures}\label{Generalized complex structures and generalized Kahler structures}\label{GCS GK}
\subsection{Generalized complex structures and nondegenerate, pure spinors}
Let $M$ be a differentiable manifold of real dimension $2n$.
The bilinear form $\lan\,\,,\,\,\ran_{\tt}$ on 
the direct sum $T_M \oplus T^*_M$ over a differentiable manifold $M$ of dim$=2n$ is defined by 
$$\lan v+\xi, u+\eta \ran_{\tt}=\frac12\(\xi(u)+\eta(v)\),\quad  u, v\in T_M, \xi, \eta\in T^*_M .$$
Let SO$(\TT)$ be the fibre bundle over $M$ with fibre SO$(2n, 2n)$ which is 
a subbundle of End$(\TT)$  preserving the bilinear form $\lan\,\,\,,\,\,\,\ran_{\tt}$ 
 An almost \complex structure $\J$ is a section of SO$(\TT)$ satisfying $\J^2=-\id.$ Then as in the case of almost complex structures, an almost \complex structure $\J$ yields the eigenspace decomposition :
\bgn{equation}\label{eigenspace decomposition}
(T_M \oplus T^*_M)^\C =\L_\J \oplus \ol \L_\J,
\end{equation} where 
$\L_\J$ is $-\sqrt{-1}$-eigenspace and  $\ol{\L}_\J$ is the complex conjugate of $\L_\J$. 
The Courant bracket of $\TT$ is defined by 
$$
 [u+\xi, v+\eta]_{\cou}=[u,v]+{\mathcal L}_u\eta-{\mathcal L}_v\xi-\frac12(di_u\eta-di_v\xi),
 $$
 where $u, v\in T_M$ and $\xi, \eta$ is $T^*_M$.
If $\L_\J$ is involutive with respect to the Courant bracket, then $\J$ is a generalized complex structure, that is, $[e_1, e_2]_{\cou}\in \Gam(\L_\J)$  for any two elements 
 $e_1=u+\xi,\,\, e_2=v+\eta\in \Gam(\L_\J)$.
%\addtocounter{section}{-1}
%\addtocounter{section}{-1}
%\setcounter{page}{3}
%\clearpage
Let $\CL(T_M \oplus T^*_M)$ be the Clifford algebra bundle which is 
a fibre bundle with fibre the Clifford algebra $\CL(2n, 2n)$ with respect to $\lan\,,\,\ran_{\tt}$ on $M$.
Then a vector $v$ acts on the space of differential forms $\oplus_{p=0}^{2n}\w^pT^*_M$ by 
the interior product $i_v$ and a $1$-form $\t$ acts on $\oplus_{p=0}^{2n}\w^pT^*_M$ by the exterior product $\t\w$, respectively.
Thus $(\TT)^\C$ acts on differential forms. (Note that by using (\ref{eigenspace decomposition}),  $\L_\J\oplus\ol\L_\J$ also acts on differential forms.)
Then the space of differential forms gives a representation of the Clifford algebra $\CL(\TT)$ which is 
the spin representation of $\CL(\TT)$. 
Thus
the spin representation of the Clifford algebra arises as the space of differential forms $$\w^\bullet T^*_M=\oplus_p\w^pT^*_M=\w^{\even}T^*_M\oplus\w^{\odd}T^*_M.$$ 
The inner product $\lan\,,\,\ran_s$ of the spin representation is given by 
$$
\lan \a, \,\,\,\b\ran_s:=(\a\w\sig\b)_{[2n]},
$$
where $(\a\w\sig\b)_{[2n]}$ is the component of degree $2n$ of $\a\w\sig\b\in\oplus_p \w^pT^*_M$ and 
$\sig$ denotes the Clifford involution which is given by 
$$
\sig\b =\bgn{cases}&+\b\qquad \deg\b \equiv 0, 1\,\,\mod 4 \\ 
&-\b\qquad \deg\b\equiv 2,3\,\,\mod 4\end{cases}
$$
We define $\ker\Phi:=\{ e\in (T_M\oplus T^*_M)^\C\, |\, e\cdot\Phi=0\, \}$ for a complex differential form $\Phi
\in \w^{\even/\odd}T^*_M.$
If $\ker\Phi$ is maximal isotropic, i.e., $\dim_\C\ker\Phi=2n$, then $\Phi$ is called {\it a pure spinor} of even/odd type.
A pure spinor $\Phi$ is {\it nondegenerate} if $\ker\Phi\cap\ol{\ker\Phi}=\{0\}$, i.e., 
$(T_M\oplus T^*_M)^\C=\ker\Phi\oplus\ol{\ker\Phi}$.
Then a nondegenerate, pure spinor $\Phi\in \w^\bullet T^*_M$ gives an almost generalized complex structure $\J_{\Phi}$ which satisfies 
$$
\J_\Phi e =
\bgn{cases}
&-\sqrt{-1}e, \quad e\in \ker\Phi\\
&+\sqrt{-1}e, \quad e\in \ol{\ker\Phi}
\end{cases}
$$
%\setcounter{page}{3}
%\clearpage
Conversely, an almost \complex structure $\J$ locally arises as $\J_\Phi$ for a nondegenerate, pure spinor $\Phi$ which is unique up to multiplication by
non-zero functions.  Thus an almost \complex structure yields the canonical line bundle $K_{\J}:=\C\lan \Phi\ran$ which is a complex line bundle locally generated by a nondegenerate, pure spinor $\Phi$ satisfying 
$\J=\J_\Phi$.
A \complex structure 
$\J_\Phi$ is integrable if and only if $d\Phi=\eta\cdot\Phi$ for a section $\eta\in (T_M\oplus T^*_M)^\C$. 
The {\it type number} of $\J=\J_\Phi$ is defined as the minimal degree of the differential form $\Phi$. Note that type number Type $\J$ is a function on a manifold which is not a constant in general.
\bgn{example}
Let $J$ be a complex structure on a manifold $M$ and $J^*$ the complex structure on the dual  bundle $T^*_M$ which is given by $J^*\xi(v)=\xi (Jv)$ for $v\in T_M$ and $\xi\in T^*_M$.
Then a \complex structure $\J_J$ is given by the following matrix
$$\J_J=\bgn{pmatrix}J&0\\0&-J^*
\end{pmatrix},$$
Then the canonical line bundle is the ordinary one which is generated by complex forms of type $(n,0)$.
Thus we have  Type $\J_J =n.$
\end{example}
\bgn{example}
Let $\ome$ be a symplectic structure on $M$ and $\h\ome$ the isomorphism from $T_M$ to $T^*_M$ given by $\h\ome(v):=i_v\ome$. We denote by $\h\ome^{-1}$ the inverse map from $T^*_M$ to $T_M$.
Then a \complex structure $\J_\psi$ is given by the following
$$\J_\psi=\bgn{pmatrix}0&-\h\ome^{-1}\\
\h\ome&0
\end{pmatrix},\quad\text{\rm Type $\J_\psi =0$}$$
Then the canonical line bundle is given by the differential form $\psi=e^{-\sqrt{-1}\ome}$. 
Thus Type $\J_\psi=0.$
\end{example}
\bgn{example}[$b$-field action]
A real $d$-closed $2$-form $b$ acts on a \complex structure by the adjoint action of Spin group $e^b$ which provides
a \complex structure $\Ad_{e^b}\J=e^b\circ \J\circ e^{-b}$. 
\end{example}
\bgn{example}[Poisson deformations]\label{Poisson deformations}
Let $\b$ be a holomorphic Poisson structure on a complex manifold. Then the adjoint action of Spin group $e^\b$ gives deformations of new \complex structures by 
$\J_{\b t}:=\Ad_{\b^{Re} t}\J_J$.  Then Type ${\J_{\b t}}_x=n-2$ (rank of $\b_x$) at $x\in M$,
which is called the Jumping phenomena of type number.
\end{example}
Let $(M, \J)$ be a generalized complex manifold and $\ol \L_\J$ the eigenspace of eigenvalue $\sqrt{-1}$.
Then we have the Lie algebroid complex $\w^\bullet\ol{\L}_\J$ (cf. \cite {Gua_2011}):
$$
0\arrow\w^0\ol \L_\J\overset{\ol\pa_\J}\arrow\w^1\ol \L_\J\overset{\ol\pa_\J}\arrow\w^2\ol \L_\J\overset{\ol\pa_\J}\arrow\w^3\ol \L_\J\arrow\cdots
$$
Since the Lie algebroid complex is an elliptic complex, 
the cohomology $H^p(\w^\bullet\ol\L_\J)$ of the Lie algebroid complex is finite dimensional if $M$ is compact.
Let $\{e_i\}_{i=1}^n$ be a local basis of $\L_\J$ for an almost \complex structure $\J$, 
where $\lan e_i, \ol e_j\ran_{\tt}=\del_{i,j}$.
The almost \complex structure $\J$ is written as an element of Clifford algebra,
$$
\J=\frac{\sqrt{-1}}2\sum_i e_i\cdot\ol {e}_i,
$$
where $\J$ acts on $\TT$ by the adjoint action $[\J, \,]$. 
Thus we have $[\J, e_i]=-\sqrt{-1}e_i$ and $[\J, \ol e_i]=\sqrt{-1}e_i$.
An almost \complex structure $\J$ acts on differential forms by the Spin representation which gives the decomposition into eigenspaces:
\bgn{equation}\label{eq:bunkai}
\w^\bullet T^*_M=U_\J^{-n}\oplus U_\J^{-n+1}\oplus\cdots U_\J^{n},
\end{equation}
where $U^{i}(=U^i_\J)$ denotes the $i$-eigenspace. Then $K_\J =U^{-n}_\J$ and 
$U_\J^{-n+p}$ is given by $\w^p\ol\L_\J \cdot K_\J$ which denotes the spin action of $\w^p\ol\L_\J$ 
on $K_\J.$ Since $\J$ is integrable, the exterior derivative $d$ is decomposed into 
$\del_\J + \ol\del_\J$, where $\del_\J: U_\J^i\to U_\J^{i-1}$ and 
$\ol\del_\J : U_\J^i\to U_\J^{i+1}.$
\subsection{Generalized K\"ahler structures}
%$M$: manifold of dim$_\R=2n$ \\
\bgn{definition}
{\it A generalized K\"ahler structure} is a pair $(\J_1, \J_2)$ consisting of two commuting \complex structures 
$\J_1, \J_2$ such that $\h G:=-\J_1\circ\J_2=-\J_2\circ \J_1$ gives a positive definite symmetric form 
$G:=\lan \h G\,\,,  \,\,\ran$ on $T_M\oplus T_M^*$, 
We call $G$ {\it a generalized metric}.
\end{definition}
\bgn{example}
Let $X=(M, J,\ome)$ be a K\"ahler manifold. Then the pair $(\J_J, \J_\psi)$ is a generalized K\"ahler where 
$\psi=\exp(\sqrt{-1}\ome)$. 
\end{example}
\bgn{example}
Let $(\J_1, \J_2)$ be a generalized K\"ahler structure. 
Then the action of $b$-fields gives a generalized K\"ahler structure 
$(\Ad_{e^b}\J_1, \Ad_{e^b}\J_2)$ for a real $d$-closed $2$-form $b.$
\end{example}
\bgn{definition}\label{Generalized Kahler structure of symplectic type}
 {\it A generalized K\"ahler structure of symplectic type} is a generalized K\"ahler structure $(\J, \J_\psi),$
where $\J_\psi$ is a generalized complex structure induced from a $d$-closed, nondegenerate, pure spinor $\psi =e^{b-\sqrt{-1}\ome}$ for a $d$-closed $2$-form $b$ and a symplectic structure $\ome.$
\end{definition}
Let $(\J_1, \J_2)$ be a generalized K\"ahler structure. Then
each $\J_i$ gives the decomposition $(\TT)^\C=\L_{\J_i}\oplus\ol \L_{\J_i}$ for $i=1,2$.
Since $\J_1$ and $\J_2$ are commutative, we have the simultaneous eigenspace decomposition 
$$
(\TT)^\C=(\L_{\J_1}\cap \L_{\J_2})\oplus (\ol \L_{\J_1}\cap \ol \L_{\J_2})\oplus (\L_{\J_1}\cap \ol \L_{\J_2})\oplus
(\ol \L_{\J_1}\cap \L_{\J_2}).
$$
Since $\h G^2=+\id$,
The generalized metric $\h G$ also gives the eigenspace decomposition: $\TT=C_+\oplus C_-$, 
where $C_\pm$ denote the eigenspaces of $\h G$ of eigenvalues $\pm1$. 
We denote by $\L_{\J_1}^\pm$ the intersection $\L_{\J_1}\cap C^\C_\pm$. 
Then it follows 
\bgn{align}\label{GK decomposition}
&\L_{\J_1}\cap \L_{\J_2}=\L_{\J_1}^+,  \quad \ol \L_{\J_1}\cap \ol \L_{\J_2}=\ol \L_{\J_1}^+\\
&\L_{\J_1}\cap \ol \L_{\J_2}=\L_{\J_1}^-,\quad \ol \L_{\J_1}\cap \L_{\J_2}=\ol \L_{\J_1}^-
\end{align}
Then $(\w^i\ol\L_{\J_1}^+)\w(\w^j\ol\L_{\J_1}^-)$ acts on $K_\J$ by the spin action to yield 
$U^{-n+i+j, i-j}:=(\w^i\ol\L_{\J_1}^+)\w(\w^j\ol\L_{\J_1}^-)\cdot K_\J.$
We have the decomposition of differential forms: 
$$
\w^\bullet T^*_M=\oplus U^{p,q}
$$
The exterior differential $d$ is also decomposed into $\del_+ +\del_- + \ol\del_+ +\ol\del_-,$
where $\del_\J =\del_+ + \del_-$ and $\ol\del_\J=\ol\del_++\ol\del_-,$
and $\del_+: U^{p,q}\to U^{p-1, q-1},\,\,\, \del_-: U^{p,q}\to U^{p-1, q+1}$ and 
$\ol\del_+: U^{p,q}\to U^{p+1, q+1}, \,\,\, \ol\del_-: U^{p,q}\to U^{p+1, q-1}.$
The generalized metric $G$ gives the formal adjoint operators $d^*$, $\del_\J^*$, \,$\ol\del_\J^*$ and $\del_\pm^*$ , 
$\ol\del_\pm^*.$ 
Then the generalized K\"ahler identity holds : $\del_+=-\del_+^*, \quad \del_-^*=\del_-.$
We denote by $\trian:=dd^*+d^*d$ the Laplacian of $d$ and $\square_{\pa_J}:=\pa_\J\pa_\J^*+\pa_\J^*\pa_\J$ the Laplacian of $\pa_\J.$
We also define the Laplacians $\square_{\ol\pa_\J}:=\ol\pa_\J\ol\pa_J^*+\ol\pa_\J^*\ol\pa_\J$
and $\square_{\del_\pm}:=\del_\pm \del_\pm^*+\del_\pm^*\del_\pm$ and 
$\square_{\ol\del_\pm}:=\ol\del_\pm \ol\del_\pm^*+\ol\del_\pm^*\ol\del_\pm.$
Then we have 
$$
\trian=2\square_{\pa_\J}=\square_{\ol\pa_\J}=4\square_{\del_\pm}=4\square_{\ol\del_\pm}
$$
Thus we have the generalized Hodge decomposition:
\bgn{proposition}[Gualtieri, \cite{Gua_2004}]
$$H^\bullet(M, \C)=H^{p,q}(M, \J_1, \J_2),$$
where $H^\bullet(M,\C)={\bigoplus_{i=0}^{\dim_\R M} H^i(M,\C)}$ and $H^{p,q}(M, \J_1, \J_2):=\ker\trian\cap U^{p,q}.$
\end{proposition}
\bgn{remark}
The decomposition does hold only when we consider cohomologies of all degrees.
\end{remark}

\subsection{The deformation-stability theorem of generalized K\"ahler manifolds}
\label{The stability theorem of generalized Kahler manifolds}
It is known that the deformation-stability theorem of ordinary K\"ahler manifolds holds
\bgn{theorem}[Kodaira-Spencer]
Let $X=(M,J)$ be a compact K\"ahler manifold and $X_t$ small deformations of $X=X_0$ as complex manifolds.
Then $X_t$ inherits a K\"ahler structure. 
\end{theorem}
The following deformation-stability theorem of generalized K\"ahler structures provides many interesting examples of generalized K\"ahler manifolds of symplectic type.
\bgn{theorem}[Goto, \cite{Goto_2010}]\label{deformation-stability theorem}
Let $X=(M,J,\ome)$ be a compact K\"ahler manifold and $(\J_J, \J_\psi)$ the induced generalized K\"ahler structure, 
where $\psi=e^{-\sqrt{-1}\ome}$. 
If there are analytic deformations $\{\J_t\}$ of $\J_0=\J_J$ as \complex structures, then there are deformations of $d$-closed nondegenerate, pure spinors $\{\psi_t\}$ such that 
pairs $(\J_t, \J_{\psi_t})$ are generalized K\"ahler structures, where $\psi_0=\psi$
\end{theorem}
Then we have the following:
\bgn{corollary}
Let $X=(M,J, \ome)$ be a compact K\"ahler manifold with a nontrivial holomorphic Poisson structure $\b.$
Then there exist  nontrivial deformations of generalized K\"ahler structures $(\J_{\b t}, \, \J_{\psi_t})$ such that 
$\{\J_{\b t}\}$ is the Poisson deformations given by Example \ref{Poisson deformations} and $\{\psi_t\}$ is a family of $d$-closed nondegenerate, pure spinors and $\psi_0=e^{-\sqrt{-1}\ome}.$
\end{corollary}
%%%%%%%%%%%%%%%%%%%%%%%%%%%%%%%%%%%%%%%%%%%%%%%%%%%%%%%

\section{Generalized Hamiltonian diffeomorphisms}\label{GHamDiff}
Let $\J$ be a generalized complex structure on a manifold $M.$
Then $\J$ acts on an exact $1$-form $du$ to give $\J du \in \TT$ for a real function $u$ on $M.$
Then we define $\frak{ham}_{\J}(M)$ by 
$$
\ham_{\J}(M):=\{\, \J du\, |\, u\in C^\infty(M,\R)\,\}
$$
%Each $\J du\in \TT$ generates an emelemt of $\wtil\Diff(M).$
The Courant bracket  on $\TT$ does not satisfies the Jacobi identity in general.
However if we restrict the Courant bracket to $\ham_\J(M),$  the Jacobi identity does hold and  we obtain a Lie algebra.
\bgn{proposition}
$\ham_{\J}(M)$ is a Lie algebra with respect to the Courant bracket.
\end{proposition}
\bgn{proof}
Since $\lan \J du_1, \, \J du_2\ran_{\tt}=\lan du_1, \, du_2\ran_{\tt}=0,$
it follows $\ham_\J(M)$ is isotropic.
Since $\J$ is integrable, the Nijenhuis tensor vanishes, 
\bgn{align}
[\J du_1, \, \J du_2]_{\cou}=&[du_1, \, du_2]_{\cou}+\J[du_1, \, \J du_2]_{\cou}+
\J[\J du_1, \, du_2]_{\cou}\\
=&\J[du_1, \, \J du_2]_{\cou}+\J[\J du_1, \, du_2]_{\cou}
\end{align}
From the definition of the Courant bracket, we have 
$$
[\J du_1,\, du_2]_{\cou}=\L_{\J du_1}( du_2) = d\L_{\J du_1}u_2
$$
We denote by $\{ u_1, u_2\}_\J$ a real function $\L_{\J du_1}u_2-\J_{\J du_2}u_1$, which reduces to the usual 
$\ome$-Poisson bracket if $\J$ is given by a nondegenerate, pure spinor $e^{-i\ome}$.
Then we obtain 
$$
[\J du_1, \, \J du_2]_{\cou}=\J d \{ u_1, \, u_2\}_\J\in \ham
_\J(M)
$$
Thus $\ham_\J(M)$ is closed under the courant bracket and isotropic. 
Hence $\ham_\J(M)$ is a Lie algebra. 
\end{proof}
Then $\ham_\J(M)$ is identified with $C^\infty(M, \R)_0:=C^\infty(M, \R)/\{\text{\rm constants }\}$.

\bgn{definition}
The Lie algebra $\ham_\J(M)$ defines 
a  connected Lie group Ham$_{\J}$ which is called 
{\it a generalized Hamiltonian diffeomorphisms} with respect to $\J.$
\end{definition}

\bgn{remark}
Let $\ome$ be a symplectic structure on $M.$ Then $e^{\sqrt{-1}\ome}$ is a $d$-closed nondegenerate, pure spinor.
If $\J$ is induced from the structure $e^{\sqrt{-1}\ome},$ then $\ham_\J(M)$ coincides with the Lie algebra of the ordinary 
Hamiltonian diffeomorphisms
\end{remark}

\section{Generalized scalar curvature as moment map}\label{GSCasmoment}
Let ${{\mathcal B}}(M)$ be the set of almost \complex structures on a differentiable compact manifold $M$ of dimension $2n$,
 that is, 
$${{\mathcal B}}(M):=\{ \J\,:\,\text{\rm almost \complex structure on }M\,\}.$$
We also define ${\mathcal B}^{\integral}(M)$ as the set of  \complex structures on $M$, i.e., integral ones
$${\mathcal B}^{\integral}(M):=\{ \J\, :\text{\rm  \complex structure on }M\,\}.$$
We fix a \complex structure $\J_\psi$ which is defined by 
a set of nondegenerate, pure spinors $\psi:=\{\psi_\a\}$ relative to a cover $\{U_\a\}$ of $M.$
Then we have 
\bgn{equation}\label{dpsia=zetaacdotpsi}
d\psi_\a =\zeta_\a\cdot\psi_\a,
\end{equation}
where we take $\zeta_\a\in \sqrt{-1}(\TT).$
We can take $\{\psi_\a\}$ which satisfies 
$$\lan \psi_\a, \,\,\ol\psi_\a\ran_s=
\lan \psi_\b, \,\,\ol\psi_\b\ran_s$$ if $U_\a\cap U_\b=\emptyset$. 
Then we define a volume form vol$_M$ to be $(\sqrt{-1})^n\lan \psi_\a, \,\,\ol\psi_\a\ran_s$
for each $\a$ which is globally defined.
An almost \complex structure $\J$ is {\it $\J_\psi$-compatible} if and only if the pair $(\J, \J_\psi)$
is an almost generalized K\"ahler structure.
Let ${\mathcal B}_{\J_\psi}(M)$ be the set of $\J_\psi$-compatible almost \complex structure, that is  
$${\mathcal B}_{\J_\psi}(M):=\{\, \J\in {{\mathcal B}(M)}\, :\, (\J,\J_\psi)\text{\rm is an almost generalized K\"ahler structure}\, \}.$$
We assume that $\B_{\J_\psi}(M)$ is not an empty set through this paper.
We also define ${\mathcal B}^{\integral}_{\J_\psi}(M)$ to be the set of $\psi$-compatible \complex structures,
For each point $x\in M$, we define ${\mathcal B}_{\J_\psi}(M)_x$ to be the set of $\psi_x$-compatible almost \complex structures on $T_xM\oplus T^*_xM$ , that is, 
 $${{\mathcal B}}_{\J_\psi}(M)_x:=\{\, \J_x\, |(\J_x, \J_{\psi, x}): \text{\rm almost generalized K\"ahler structure at } x \, \}.$$
 Then we see that ${\mathcal B}_{\J_\psi}(M)_x$ is given by the Riemannian Symmetric space of type  AIII
$$U(n,n)/U(n)\times U(n)$$ which is biholomorphic to the complex bounded domain 
 $$\{\, h\in M_n(\C)\, |\, 1_n-h^*h>0\, \},$$ where $M_n(\C)$ denotes the set of complex matrices of $n\times n.$ 
 \bgn{remark} In K\"ahler geometry, the set of almost complex structures compatible with a symplectic structure $\ome$ is given by 
 the Riemannian symmetric space Sp$(2n)/U(n)$ which is biholomorphic to the Siegel upper half plane 
 $$\{\, h\in \GL_n(\C)\, |\, 1_n-h^*h>0,\, h^t=h \,\}$$
 \end{remark}
Let $P_{\J_\psi}$ be the fibre bundle over $M$ with fibre ${{\mathcal B}}_{\J_\psi}(M)_x$, that is, 
$$P_{\J_\psi}:=\bigcup_{x\in M}{{\mathcal B}_{\J_\psi}(M)_x}\to M,$$
Then ${\mathcal B}_{\J_\psi}(M)$ is given by smooth sections $\Gam (M, P_{\J_\psi})$ which contains the integral ones ${\mathcal B}^{\integral}_{\J_\psi}(M)$. 
We can introduce a Sobolev norm on ${\mathcal B}_{\J_\psi}(M)$ such that ${\mathcal B}_{\J_\psi}(M)$ becomes a Banach manifold in the standard method.
The tangent bundle of ${\mathcal B}_{\J_\psi}(M)$ at $\J$ is given by 
$$T_{\J}{\mathcal B}_{\J_\psi}(M)=\{\, \dot{\J}\in\text{\rm so}(T_M\oplus T^*_M)\,:\, \dot{\J}\J+\J\dot{\J}=0,\, \dot{\J}\J_\psi=\J_\psi\dot{\J}\, \},$$
where so$(\TT)$ denotes the set of sections of Lie algebra bundle of SO$(\TT)$.
Then it follows that there exists an almost complex structure $J_{\mathcal B}$ on 
${\mathcal B}_{\J_\psi}(M)$which is given by 
$$
J_{\mathcal B}(\dot{\J}):=\J\dot{\J}, \qquad (\,\,\dot{\J}\in T_{\J}{\mathcal B}_{\J_\psi}(M) \,\,)
$$
We also have a Riemannian metric $g_{\mathcal B}$ and a $2$-form $\Ome_{\mathcal B}$ on 
${\mathcal B}_{\J_\psi}(M)$ by 
\bgn{align}\label{Apsi}
&g_{\mathcal B}(\dot{\J_1},\dot{\J_2}):=\int_M \tr(\dot{\J_1}\dot{\J_2})\,\vol_M
\\
&\Ome_{\mathcal B}(\dot{\J_1},\dot{\J_2}):=-\int_M \tr(\J\dot{\J_1}\dot{\J_2})
\,\vol_M
\end{align}
for $\dot{\J_1}, \dot{\J_2}\in T_{\J}{\mathcal B}_{\J_\psi}(M)$. 
\bgn{proposition}
$J_{\mathcal B}$ is integrable almost complex structure on ${\mathcal B}_{\J_\psi}(M)$ and
$\Ome_{{\mathcal B}}$ is a K\"ahler form on ${\mathcal B}_{\J_\psi}(M).$
\end{proposition}
\bgn{proof}
Let $\J_{V}$ be an almost generalized complex structure on a real vector space $V$ of dimension $2n$.
We denote by $X_n$ the Riemannian symmetric space  
$U(n,n)/U(n)\times U(n)$ which is identified with the set of almost generalized complex structures
compatible with $\J_{V}.$
We already see that ${\mathcal B}_{\J_\psi}(M)$ is the set of global sections of the fibre bundle $P_{\J_\psi}$ over a manifold $M$ with fibre 
$X_n$ which is biholomorphic to the bounded domain $\{\, h\in M_n(\C)\, |\, 1_n-h^*h>0\, \}.$
Let $\J_0$ be an element of $\B_{\J_\psi}(M)$. 
Then a generalized K\"ahler structure $(\J_0, \J_\psi)$ gives the decomposition of $(\TT)^\C$ as in (\ref{GK decomposition})
$$
(\TT)^\C=\L_{\J_\psi}^+\oplus \L_{\J_\psi}^-\oplus\ol\L_{\J_\psi}^+\oplus \ol\L_{\J_\psi}^-,
$$
where $\L_{\J_\psi}^+=\L_{\J_\psi}\cap\L_{\J_0}$ and $\L_{\J_\psi}^-=\ol\L_{\J_\psi}\cap \L_{\J_0}$.
Note that the adjoint action of the  group $ SO(\TT)$ on  the set of almost generalized complex structures is transitive.
An element of SO$(\TT)$ preserves $\J_\psi$ if and only if it preserves $\L_{\J_\psi}$.
Then every real element $g$ of SO$(\TT)$ preserving $\J_\psi$ is given by the following 
$$
g=
\bgn{pmatrix}
A&C\\B&D
\end{pmatrix}\in \text{\rm GL}(\L_{\psi}),
$$
where $A\in $End$(\L_{\J_\psi}^+)$ and  $B\in $Hom$(\L_{\J_\psi}^+, \L_{\J_\psi}^-)$ 
and $C\in$Hom$(\L_{\J_\psi}^-, \L_{\J_\psi}^+)$ and $D\in $\End$(\L_{\J_\psi}^-)$ satisfy 
\bgn{align}
&A^*A-B^*B=\id_{\L_{\J_\psi}^+}, \quad -C^*C+D^*D=-\id_{\L_{\J_\psi}^-}\\
&A^*D=B^*C
\end{align}
Thus it follows that $g$ is a section of the fibre bundle whose fibre is U$(n,n)$.
Let $h:=BA^{-1}\in $Hom$(\L_{\J_\psi}^+, \, \L_{\J_\psi}^-)$. Then $h$ satisfies 
\bgn{align}\label{idLpsi+-hh*>0}
\id_{\L_{\J_\psi}^+}-h^*h>0\in \text{End}(\L_{\J_\psi}^+)
\end{align}
Thus it follows that the fibre bundle $P_{\J_\psi}$ is identified with an open fibre bundle with fibre the bounded domain, which is a open subbundle of the complex vector bundle Hom$(\L_{\J_\psi}^+,\,\L_{\L_\psi}^-)$.
Thus $\B_{\J_\psi}(M)$ is the space of global sections of an fibre bundle over 
 $M$ which is an open subbundle of a complex vector bundle Hom$(\L_{\J_\psi}^+,\,\L_{\L_\psi}^-)$ over $M$.
In general, the set of global sections of the complex vector bundle is a complex manifold, cf. \cite{Pa_2016}. (We choose a Sobolev norm $L^2_k$ to consider the set of $L^2_k$-sections.)
Since ${\mathcal B}_{\J_\psi}(M)$ is an open set of the global sections of the complex vector bundle, 
${\mathcal B}_{\J_\psi}(M)$ is a complex submanifold with a complex structure $J'_{\mathcal B}$.
Since the almost complex structure $J_{\mathcal B}$ is induced from the complex structure of the complex bounded domain, it follows that $J_{\mathcal B}$ coincides with $\J'_{\mathcal B}$. 
Thus it follows $\J_{\mathcal B}$ is integrable almost complex structure on $\mathcal B_{\J_\psi}(M)$.
We denote by $g_{\scriptscriptstyle X_n}$ the Riemannian metric on $X_n$ and by $\ome_{\scriptscriptstyle X_n}$ the K\"ahler form which are respectively given by 
$$
g_{\scriptscriptstyle X_n}(\dot\J_1, \dot\J_2)=\tr(\dot\J_1\dot\J_2)
$$
$$
\ome_{\scriptscriptstyle X_n}(\dot\J_1, \dot\J_2)=-\tr(\J \dot\J_1\dot\J_2),
$$
where $\dot\J_1, \dot\J_2\in T_{\J}X_n$. The complex bounded domain $\{\, h\in \GL_n(\C)\, |\, 1_n-h^*h>0\, \}$ admits a K\"ahler structure which is given 
by 
$$
4\sqrt{-1}\pa\ol\pa \log\det (1_n-h^*h).
$$
Then under the identification $X_n\cong \{\, h\in M_n(\C)\, |\, 1_n-h^*h>0\, \}$ by using $\J_V,$ we have 
$\ome_{X_n}=4\sqrt{-1}\pa\ol\pa \log\det (1_n-h^*h).$
Since $\mathcal B_{\J_\psi}(M)$ is the set of global section of the open subbundle, 
the tangent bundle $T\mathcal B_{\J_\psi}(M)$ is canonically identified with the trivial bundle 
$\mathcal B_{\J_\psi}(M)\times \Gam(M, \Hom(\L_{\J_\psi}^+, \, \L_{\J_\psi}^-)).$
The complex manifold ${\mathcal B}_{\J_\psi}(M)$ inherits a Riemannian metric $g_{\mathcal B}$ and a Hermitian $2$-form $\Ome_{\mathcal B}$ which are given by 
\bgn{align}
g(\dot\J_1 , \,\,\dot\J_2):=&\int_M \tr(\dot\J_1\dot\J_2)\,\vol_M
\\ 
\Ome_{\mathcal B}(\dot\J_1 , \,\,\dot\J_2):=&-\int_M \tr(\J \dot\J_1\dot\J_2)\,\vol_M \label{OmemathcalB}
\end{align}
Since the tangent bundle of $\B_{\J_\psi}(M)$ is canonically identified with the trivial bundle, 
each global section $a\in \Hom(\L_{\J_\psi}^+, \, \L_{\J_\psi}^-)$ gives a vector field $\dot\J_\a$ of $\B_{\J_\psi}(M)$ such that $\dot\J(\J)=a$ for all $\J\in \B_{\J_\psi}(M)$.
By using $L_k^2$-metric, we obtain a basis of vector fields of $\B_{\J_\psi}(M)$ by using global sections of $\Hom(\L_{\J_\psi}^+, \, \L_{\J_\psi}^-)$. 
We also denote by $\dot\J_{a, x}$ the vector field on the fibre at $x\in M$.
For global sections $a_1, a_2\in \Hom(\L_{\J_\psi}^+, \, \L_{\J_\psi}^-)$ and for each $x\in M$, 
we have 
$${4}\sqrt{-1}\,\pa\ol\pa\log\det (1_n-h^*h) (\dot\J_{a_1, x}, \,\,\dot\J_{a_2, x})=-\tr(\J \dot\J_{a_1, x}\dot\J_{a_2, x}).$$
Then it follows 
\bgn{align}\label{OmeBdotJa1dotJa2)}
\Ome_\B(\dot\J_{a_1}, \,\,\dot\J_{a_2})
=&
 \int_M{4}\sqrt{-1}\,\pa\ol\pa\log\det (1_n-h^*h) (\dot\J_{a_1}, \,\,\dot\J_{a_2})\,\vol_M
\end{align}
Since $\B_{\J_\psi}(M)$ is given by the set of global sections of $\Hom(\L_{\J_\psi}^+, \, \L_{\J_\psi}^-)$ satisfying 
(\ref{idLpsi+-hh*>0}), 
$h\mapsto \int_M \log\det (\id_{\L_{\J_\psi}^+}-h^*h)\vol_M$ is regarded as a function on $\B_{\J_\psi}(M)$. 
Let $\ol\pa_{\B}$ be the $\ol\pa$-operator of the complex manifold $\B_{\J_\psi}(M)$ and 
$\pa_{\B}$ the complex conjugate of $\ol\pa_{\B}$. 
Then 
$$
\(\pa_\B\ol\pa_\B \int_M \log\det (\id_{\L_{\J_\psi}^+}-h^*h)\vol_M\) (\dot\J_{a_1}, \,\, \dot\J_{a_2})
$$ is given by $(1,1)$-component of 
$$\frac{d}{dt_1}\frac{d}{dt_2}\big|_{t_1, t_2=0}
\int_M  \log\det (\id_{\L_{\J_\psi}^+}-h_{t_1, t_2}^*h_{t_1, t_2})\vol_M,
$$
where $h_{t_1, t_2}=h+a_1 t_1+ a_2t_2\in \Gam(M, \Hom(\L_{\J_\psi}^+, \,\L_{\J_\psi}^-)).$ and $t_1, t_2$ are parameters of small deformations.
Then it follows 
\bgn{align}
\frac{d}{dt_1}\frac{d}{dt_2}\big|_{t_1, t_2=0}
\int_M  \log\det (\id_{\L_{\J_\psi}^+}-h_{t_1, t_2}^*h_{t_1, t_2})\vol_M\\
=
\int_M \frac{d}{dt_1}\frac{d}{dt_2}\big|_{t_1, t_2=0} \log\det (\id_{\L_{\J_\psi}^+}-h_{t_1, t_2}^*h_{t_1, t_2})\vol_M,
\end{align}
Since $\pa\ol\pa\log\det (1_n-h^*h) (\dot\J_{a_1}, \,\,\dot\J_{\a_2})$ is given by the $(1,1)$-component of 
$$\frac{d}{dt_1}\frac{d}{dt_2}\big|_{t_1, t_2=0} \log\det (\id_{\L_{\J_\psi}^+}-h_{t_1, t_2}^*h_{t_1, t_2}).$$
From (\ref{OmeBdotJa1dotJa2)}), we have 
$$
\Ome_\B(\dot\J_{a_1}, \,\,\dot\J_{a_2})=
4\sqrt{-1}\(\pa_\B\ol\pa_\B \int_M \log\det (\id_{\L_{\J_\psi}^+}-h^*h)\vol_M\) (\dot\J_{a_1}, \,\, \dot\J_{a_2})
$$
Thus $\Ome_\B$ is $\pa_\B\ol\pa_\B$-exact.
Hence $\Ome_{\mathcal B}$ is closed.
Thus $({\mathcal B}_{\J_\psi}(M),J_{\mathcal B}, \,\Ome_{\mathcal B}) $ is a K\"ahler manifold.
\end{proof}
Let $\wtil{\Diff}(M)$ be an extension of diffeomorphisms of $M$ by $2$-forms which is defined as 
$$
\wtil{\Diff}(M):=\{\, e^b F\,:\, F\in \Diff(M),\,\, b: 2\text{\rm -form}\, \,\}.
$$
Note that the product of $\wtil{\Diff}(M)$ is given by 
$$
(e^{b_1}F_1)( e^{b_2}F_2) :=e^{b_1+F_1^*(b_2)}F_1\circ F_2,
$$
where $F_1, F_2\in \Diff(M)$ and $b_1, b_2$ are real $2$-forms.
The action of $\wtil{\Diff}(M)$ on ${{\mathcal GC}}(M)$ by 
\bgn{equation}\label{ebFcircJcircF-1}
e^{b} F_\#\circ \J\circ F_\#^{-1} e^{-b}, 
\end{equation}
 where $F\in \Diff(M)$ acts on $\J$ by $F_\#\circ \J\circ F_\#^{-1}$ and 
 and $e^b$ is regarded as an element of SO$(\TT)$ and $F_\#$ denotes the bundle map of $\TT$ which is the lift of $F.$
 \bgn{remark}
An element $v+\t\in \TT$ generates a $1$-parameter family of $\wtil\Diff(M)$. 
An element of the Lie algebra of $\wtil\Diff(M)$ is a 
pair $(v, d\t)$ which consists of  a vector field $v$ and a $d$-exact $2$-form $d\t.$
Then the Lie bracket is given by 
$$
[v_1+d\t_1, v_2+d\t_2]=[v_1, v_2]+\L_{v_1}\t_2-\L_{v_2}\t_1
$$
Since $\ham_J(M)$ is an isotropic subspace of 
$\TT,$ we have a homomorphism from the Lie algebras
$\ham_J(M)$ to the Lie algebra of $\wtil\Diff(M).$
\end{remark}
For a (integral) generalized complex structure $\J_\psi,$
we define $\wtil{\Diff(M)}_{\J_\psi}$ to be a subgroup consists  of elements of $\wtil{\Diff(M)}$ which preserves $\J_\psi$, 
$$
\wtil{\Diff}_{\J_\psi}(M)=\{\, e^bF\in \wtil{\Diff}(M)\, : 
e^{b} F_\#\circ \J_\psi\circ F_\#^{-1} e^{-b}=\J_\psi\, \}.
$$
Then from (\ref{Apsi}), we have the following,
\bgn{proposition}\label{OmeBinvawtilDiffJpsiM}
The symplectic structure $\Ome_{\mathcal B}$ 
is invariant under
the action of $\psi$-preserving group $\wtil{\Diff}_{\J_\psi}(M)$.\end{proposition}
\bgn{proof}
The result follows from (\ref{OmemathcalB}) and (\ref{ebFcircJcircF-1}) since $\vol_M$ is invariant under the action of $\wtil{\Diff}_{\J_\psi}(M)$.
\end{proof}
\bgn{proposition}\label{GJpsi preserves Ome}
Let $\text{\rm Ham}_{\J_\psi}$be the generalized Hamiltonian diffeomorphisms whose Lie algebra is $\ham_{\J_\psi}(M).$
Then $\text{\rm Ham}_{\J_\psi}$also preserves $\Ome_{\mathcal B}.$
\end{proposition}
\bgn{proof}
The Lie algebra $\ham_{\J_\psi}$ of the Lie group Ham$_{\J_\psi}$is given by 
$\{e:=\J_\psi du\, |\, u \in C^\infty(M,\R)\,\}$ as before. 
Then the action of $\J_\psi du$ on $\psi_\a$ is given by the Lie derivative 
$\L_{e}\psi_\a= d (\J_\psi du)\cdot\psi_\a+ (\J_\psi du)d\psi_\a.$
Since $\sqrt{-1}du +\J_\psi(du)\in \L_{\J_\psi},$ we have $(\sqrt{-1}du +\J_\psi(du))\cdot\psi_\a=0$.
From (\ref{dpsia=zetaacdotpsi}), we have
\bgn{align*}
\L_{e}\psi_\a=&-\sqrt{-1}d ((du)\psi_\a)+(\J_\psi du)\cdot\zeta_\a\cdot\psi_\a\\
=&\sqrt{-1}(du)\w d\psi_\a+(\J_\psi du)\cdot\zeta_\a\cdot\psi_\a\\
=&\sqrt{-1}(du)\w \zeta_\a\cdot\psi_\a+(\J_\psi du)\cdot\zeta_\a\cdot\psi_\a\\
=&\(\sqrt{-1}(du)+(\J_\psi du)\)\cdot\zeta_\a\cdot\psi_\a
\end{align*}
Since $\(\sqrt{-1}(du)+(\J_\psi du)\)\in \L_{\J_\psi}$, we see that the component 
$$\pi_{U_{\J_\psi}}^{-n+2}\(\sqrt{-1}(du)+(\J_\psi du)\)\cdot\zeta_\a\cdot\psi_\a=0.$$
Thus $\L_e\psi_\a$ is in $K_{\J_\psi}.$
Hence $\L_e$ preserves the canonical line bundle $K_{\J_\psi}$ and then it follows that 
$\ham_{\J_\psi}(M)$ preserves $\J_\psi.$
Thus $\text{\rm Ham}_{\J_\psi}$ also preserves $\J_\psi$ and $\vol_M$.
The infinitesimal action of $(v, \t)\in \ham_{\J_\psi}(M)$ on $\B_{\J_\psi}(M)$ is given by 
$L_v+ d\t$ which is the infinitesimal action of $(v, d\t)$ of the Lie algebra of $\wtil\Diff_{\J_\psi}(M)$.
From Proposition \ref{OmeBinvawtilDiffJpsiM}, it follows that $\Ome_\B$ vanishes by the infinitesimal action of $(v, \t)\in \ham_{\J_\psi}(M)$. 
Thus one see that the action of Ham$_{\J_\psi}$ preserves $\Ome_\B$.
\end{proof}
As is shown before,
the Lie algebra $\ham_{\J_\psi}(M)$ is given by $C_0^\infty(M),$
where $C^\infty_0(M)=\{\, f\in C^\infty(M)\, |\, \int_M f\,\vol_M=0\,\}.$
Then $e:=\J_\psi(df)\in \TT$ is called {\it a generalized Hamiltonian element}.
Note that  we have $e\cdot\psi_\a=-\sqrt{-1}df\cdot\psi_\a$.

%\clearpage
In order to show the existence of the moment map, we shall restrict our attention to generalized K\"ahler manifolds $(\J, \J_\psi)$, where 
$\J_\psi$ is induced from a set of locally defined $d$-closed nondegenerate, pure spinors $\psi:=\{\psi_\a\}$.
\bgn{theorem}\label{existence of moment map}
We assume that $\J_\psi$ is induced from a set of $d$-closed, nondegenerate, pure spinors $\psi :=\{\psi_\a\}.$
Then there exists a moment map  
$\mu: {\mathcal B}_{\J_\psi}(M)\to C^\infty_0(M,\R)^*$ for the action of 
the generalized Hamiltonian diffeomorphisms $\Ham_{\J_\psi}$, which is explicitly written in terms of pure spinors.
\end{theorem}

\bgn{remark}
In the previous paper \cite{Goto_2016}, the existence of the moment map was shown in the rather restricted cases of generalized K\"ahler manifolds of symplectic type. Thus our theorem is a generalization of the previous one and
the method of our proof is improved. 
\end{remark}
\bgn{remark}
Boulanger also obtained the moment map in the cases of toric generalized K\"ahler manifolds of symplectic type by using 
a description of toric geometry \cite{Bou_2019}. 
As though Boulanger's description of the moment map seems to be different from the one in  \cite{Goto_2016}, 
these should match each other since the moment map is unique modulo constant.
In fact, 
Yicao Wang actually shows that these are the same by using explicit calculations \cite{WY_2020}.
\end{remark}
\bgn{remark}
J.~Streets studies problems of generalized K\"ahler structures by using  pluriclosed flows \cite{St_2016}. 
A generalized K\"ahler structure of type $(0,0)$ is a generalized K\"ahler structure $(\J_1, \J_2)$ consisting of 
two generalized complex structures of type $0$,
which is also called a degenerate generalized K\"ahler structure.
In the cases of generalized K\"ahler structures of type $(0,0)$, 
his definition of 
generalized K\"ahler structure with constant scalar curvature is the same as the one in our paper 
(see also \cite{Goto_2016}, for generalized K\"ahler structures of type $(0,0)$).
The Calabi-Yau type problem of generalized K\"ahler manifolds of type $(0,0)$ was discussed by Apostolov and Streets in
\cite{AS_2017}.
\end{remark}
In order to show Theorem \ref{existence of moment map}, we need several Lemmata.
Let $\J\in\mathcal B_{\J_\psi}(M)$ be an almost generalized complex structure which is induced from 
a set of nondegenerate, pure spinors $\phi=\{\phi_\a\}.$
We normalize $\{\phi_\a\}$ such that $\lan \phi_\a,\,\,\ol\phi_\a\ran_s =\vol_M$ for each $\a.$
Then one has
\bgn{lemma}
 $d\phi_\a$ is given by  
\bgn{equation}\label{dphia=eta+N}
d\phi_\a=(\eta_\a+N_\a)\cdot\phi_\a,
\end{equation}
where  $\eta_\a\in \sqrt{-1}(\TT)$ and $N_\a\in (\w^3\L_\J\oplus \w^3\ol\L_\J)^\R.$
Note that $\eta_\a$ and $N_\a$ are uniquely determined
\end{lemma}
\bgn{proof}
It suffices to show that $d\phi_\a\in U_\J^{-n+1}\oplus U^{-n+3}_\J$.
In fact, one has 
$$
e_1\cdot e_2\cdot e_3\cdot e_4\cdot d\phi_\a=e_1\cdot e_2\cdot [e_3, e_4]_{co}\cdot\phi_\a.
$$
for any $e_1, e_2, e_3, e_4\in \L_\J$. 
Since $[e_3, e_4]_{co}$ is given by $e_5+\ol e_6$, for some $e_5\in \L_\J$ and $\ol e_6\in \ol\L_\J$, it follows from $e_2\cdot\ol e_6+\ol e_6\cdot e_2=2\lan e_2, \,\ol e_6\ran_{\tt}$ and and 
$\L_\J=\ker \phi_\a$ that 
$$
e_1\cdot e_2\cdot e_3\cdot e_4\cdot d\phi_\a=e_1\cdot e_2\cdot \ol e_6\cdot\phi_\a=
2\lan e_2,\,\ol e_6\ran_{\tt} \, e_1\cdot\phi_\a=0.
$$
Thus one has $d\phi_\a= \eta'_\a \cdot\phi_\a+N'_\a\cdot\phi_\a \in U_\J^{-n+1}\oplus U^{-n+3}_\J$, 
where $\eta'_\a\in \ol\L_\J$ and $N'_\a\in \w^3\ol\L_\J$.
Then $\eta_\a$ is the imaginary element $\eta'_\a-\ol{\eta'}_\a$ and $N_\a$ is the real one $N'_\a+\ol{N'}_\a$.
\end{proof}
\bgn{remark} Note that $N_\a$ is a real element and  $N_\a=N_\b$ for all $\a, \b$. 
Then $N_\a$ defines a global element $N$, which is called {\it Nijenhuis tensor.}
In fact, $\J$ is integrable if and only if $N$ vanishes.
\end{remark}
\bgn{lemma}\label{Ncdotpsia=0}
$N\cdot\psi_\a=0$ \end{lemma}
\bgn{proof}
Since $N$ is uniquely defined by (\ref{dphia=eta+N}), 
for $e_1, e_2, e_3\in \L_\J$, we have
\bgn{align}
N(e_1, e_2, e_3)\lan \phi_\a, \,\,\ol\phi_\a\ran_s=&
\lan d\phi_\a,\,\, e_1\cdot e_2\cdot e_3\cdot\ol\phi_\a\ran_s=
-\lan e_1\cdot e_2\cdot d\phi_\a,\,\, e_3\cdot\ol\phi_\a\ran_s\\
=&-\lan [e_1, e_2]_{\cou}\cdot\phi_\a,\,\,e_3\cdot\ol\phi_\a\ran_s
\end{align}
Thus we have 
\bgn{equation}\label{N(e1,e2,e3)}
N(e_1, e_2, e_3)=2\lan [e_1, e_2]_{\cou}, \,\, e_3\ran_{\tt}
\end{equation}
This implies that $N=0$ if and only if $\J$ is integrable. 
By using $\J_\psi,$ we have the decomposition 
$\L_\J =\L_\J^+\oplus \L_\J^-$ and $\ol\L_\J =\ol\L_\J^+\oplus\ol\L_\J^-.$
Since $\ker \psi_\a =\L_\J^+\oplus\ol\L_\J^-$ and $N\in (\w^3\L_\J\oplus \w^3\ol\L_\J)^{\Re},$ we have 
$N\cdot\psi =(\ol N^{+}+ N^-)\cdot\psi,$
where $\ol N^+\in \w^3\ol\L_\J^+$ and $N^-\in \w^3\L_\J^-.$
From (\ref{N(e1,e2,e3)}), we see 
$$
N(e_1^-, e_2^-, e_3^-)=\lan [e_1^-, e_2^-]_{\cou},\,\, e_3^-\ran_{\tt},\qquad\( e_1^-, e_2^-, e_3^-\in \ol\L_{\J}^-\).
$$
Since $\J_\psi$ is integrable, it follows that 
$[e_1^-, e_2^-]_{\cou}\in \L_{\J_\psi}$. Since $e_3^-\in \L_{\J_\psi},$ we have 
$N(e_1^-, e_2^-, e_3^-)=0.$ Then it follows $N^-=0$. 
We also have $\ol N^+=0.$
Hence $N\cdot\psi =0.$
\end{proof}
\bgn{lemma}\label{dotN(t)cdotpsi=0} 
Let $\J_t$ be deformations of $\J$ such that $(\J_t, \J_\psi)$ is an almost generalized K\"ahler structures. 
We denote by  $\{\phi_\a(t)\}$ a family of nondegenerate, pure spinors which gives $\J_t$ and 
$d\phi_\a(t)=(\eta_\a(t)+N(t))\cdot\phi_\a(t),$
where $\eta_\a(t)\in \sqrt{-1}(\TT)$ and $N \in (\w^3\L_{\J_t}\oplus \w^3\ol\L_{\J_t})^{\Re}.$ Let $\dot N=\frac{d}{dt}N(t)|_{t=0}$. Then we have 
$$
\dot N\cdot \psi_\a=0
$$
\end{lemma}
\bgn{proof}
From Lemma \ref{Ncdotpsia=0}, we have $N(t)\cdot\psi_\a=0$ for all $t.$
Then we have the result.
\end{proof}
\bgn{lemma}\label{lanecdotphia, dot Ncdotolphiaran=0}
$\lan e\cdot\phi_\a, \,\,\dot N\cdot\ol\phi_\a\ran_s=0.$
\end{lemma}
\bgn{proof}
The space $\w^4(\TT)$ is decomposed into $\w^4 T_M \oplus (\w^3 T_M \otimes T_M^*)
\oplus (\w^2T_M\otimes \w^2T_M^*)\oplus (T_M\otimes \w^3T_M^*)\oplus \w^4T_M^*.$
We denote by Cont$^{2,2}$ the contraction of the component $(\w^2T_M\otimes \w^2T_M^*)$
which yields a map from $\w^4(\TT)$ to $C^\infty(M).$
Then it follows 
\bgn{align}
\lan e\cdot\phi_\a, \,\,\dot N\cdot\ol\phi_\a\ran_s=&-\lan \phi_\a, \,\,e\cdot \dot N\cdot\ol\phi_\a\ran_s\\
=&-\text{\rm Cont}^{2,2}
(e\cdot \dot N)\lan\phi_\a, \,\,\ol\phi_\a\ran_s
\end{align}
Since $\lan \phi_\a, \,\,\ol\phi_\a\ran_s=\lan \psi_\a,\,\,\ol\psi_\a\ran_s,$
we have 
\bgn{align}
\lan e\cdot\phi_\a, \,\,\dot N\cdot\ol\phi_\a\ran_s=&-\text{\rm Cont}^{2,2}
(e\cdot \dot N)\lan\psi_\a, \,\,\ol\psi_\a\ran_s\\
=&\lan e\cdot\psi_\a, \,\, \dot N\cdot \ol\psi_\a\ran_s
\end{align}
Since $\dot N$ is real, it follows from Lemma \ref{dotN(t)cdotpsi=0} that 
$\dot N\cdot\ol\psi_\a=0.$
Hence we have
$\lan e\cdot\phi_\a, \,\,\dot N\cdot\ol\phi_\a\ran_s=0.$
\end{proof}

\bgn{proof} of Theorem \ref{existence of moment map}.
Every infinitesimal deformation  of $\J$ is written by the adjoint action of $h$ 
$$
\dot\J_h :=[h,\,\J]
$$
where $h$ denotes a real element $ (\w^2\L_{\J}\oplus\w^2\ol\L_{\J})^\R.$
Then the corresponding infinitesimal deformation of $\phi$ is given by the Clifford action of $h$ on each $\phi_\a,$
$$
\dot{\phi_\a}=h\cdot\phi_\a
$$
An Hamiltonian element $e=\J_\psi df$ gives the infinitesimal deformation $\L_e\J$ of $\J.$
Then the corresponding infinitesimal deformation of $\phi$ is given by 
$$
\L_e\phi_\a=d (e\cdot\phi_\a)+e\cdot d\phi_\a=d (e\cdot\phi_\a)+e\cdot (\eta_\a+N)\cdot\phi_\a.
$$
(For simplicity, $d (e\cdot\phi_\a)$ is also denoted by $de\cdot\phi_\a$. ) 
Then in order to show the existence of a moment map, we shall calculate 
$\Ome_\B(\L_e\J, \,\, \dot\J_h)$.
%%%%%%%%%%%%%%%%%%%%%% revision %%%%%%%%%%%
In \cite{Goto_2016}, the following description of $\Ome_\B$ in terms of pure spinors is already given. (see Lemma 7.1 and (7.1) in Section 7 of \cite{Goto_2016}, where $\vol_M=i^{-n}\lan \psi, \ol\psi\ran_s$ and $\rho_\a=1$.)
%%%%%%%%%%%%%%%%%%%%%%%%
\bgn{equation}\label{OmeB(dotJh1,dotJh2)}\Ome_\B(\dot\J_{h_1}, \,\,\dot\J_{h_2})=
 c_n \Im\( i^{-n}\int_M \lan h_1\cdot\phi_\a, \,\,h_2\cdot\ol\phi_\a \ran_s\),
 \end{equation}
 where $c_n$ is a constant depending only on $n$.
Applying (\ref{OmeB(dotJh1,dotJh2)}), we obtain 
\bgn{align}
\Ome_\B(\L_e\J, \,\, \dot\J_h)=&c_n \Im\(i^{-n}\int_M \lan\L_e\phi_\a, \,\, h\cdot\ol\phi_\a\ran_s\)\\
=&c_n\Im\(i^{-n}\int_M\lan d e\cdot\phi_\a+e\cdot (\eta_\a+N)\cdot\phi_\a, \,\, h\cdot\ol\phi_\a\ran_s
\)
\end{align}
Since $h\in(\w^2\L_{\J}\oplus\w^2\ol\L_\J)^\R, $
we have 
$\lan \phi_\a, \,\, h\cdot\ol\phi_\a\ran=0.$
Since $e\cdot\eta_\a+\eta_\a\cdot e =2\lan e, \,\,\eta_\a\ran_{\tt},$ we have  
$$
\lan e\cdot\eta_\a\cdot\phi_\a,\,\, h\cdot\ol\phi_\a\ran_s =-\lan \eta_\a\cdot e\cdot\phi_\a, \,\,h\cdot\ol\phi_\a\ran_s
$$
Then we have 
\bgn{align}
\Ome_\B(\L_e\J, \,\, \dot\J_h)=&c_n \Im\(i^{-n}\int_M\lan (d-\eta_\a+N) e\cdot\phi_\a, \,\, h\cdot\ol\phi_\a\ran_s\)
\end{align}
%%%%%%%%%%%%%%%%%%%%%%%%%%%%%%
It follows  $d\lan (e\cdot\phi_\a),\,\,h\cdot\ol\phi_\a\ran_{[2n-1]}=\lan d(e\cdot\phi_\a),\,\,h\cdot\ol\phi_\a\ran_s
-\lan (e\cdot\phi_\a),\,\,d(h\cdot\ol\phi_\a)\ran_s$, where  
$\lan (e\cdot\phi_\a),\,\,h\cdot\ol\phi_\a\ran_{[2n-1]}$ is the component of $(2n-1)$-form of 
$(e\cdot\phi_\a)\w\sig(h\cdot\ol\phi_\a)$.
Recall that we assume the normalization  $i^{-n}\lan \phi_\a, \,\ol\phi_\a\ran_s=\vol_M$ for all $\a$.
Thus it follows $\phi_\a=e^{ip_{\a, \b}}\phi_\b$, where $p_{\a,\b}$ denotes a real function.
Then it follows 
$(e\cdot\phi_\a)\w\sig(h\cdot\ol\phi_\a)=(e\cdot\phi_\b)\w\sig(h\cdot\ol\phi_\b)$. 
Thus $(e\cdot\phi_\a)\w\sig(h\cdot\ol\phi_\a)$ gives a globally defined $(2n-1)$-form.
%%%%%%%%%%%%%%%%%
From the Stokes Theorem, we have 
$$
\int_M \lan d(e\cdot\phi_\a),\,\,h\cdot\ol\phi_\a\ran_s=\int_M \lan e\cdot\phi_\a, \,\,d(h\cdot\ol\phi_\a)\ran_s
$$
Since $\eta_\a$ is in $\sqrt{-1}(\TT)$, we also have 
$$
\lan \eta_\a\cdot e\cdot \phi_\a, \,\, h\cdot\ol\phi_\a\ran_s=-\lan e\cdot\phi_\a, \,\, \eta_\a\cdot h\cdot\ol\phi_\a\ran_s.
$$
Since $N$ is $\w^3(\TT)^\R$, it follows 
$$
\lan N\cdot e\cdot \phi_\a, \,\, h\cdot\ol\phi_\a\ran_s=+\lan e\cdot\phi_\a, \,\, N\cdot h\cdot\ol\phi_\a\ran_s.
$$
Substituting them, we obtain 
$$
\Ome_\B(\L_e\J, \,\, \dot\J_h)=c_n\Im\(i^{-n}\int_M\lan e\cdot\phi_\a, \,\, (d+\eta_\a-N)\cdot (h\cdot\ol\phi_\a)\ran_s\)
$$
Let $\phi(t)=\{\phi_\a(t)\}$ be a family of nondegenerate, pure spinors which gives 
$$
\frac{d}{dt}\phi_\a(t)|_{t=0}=h\cdot\phi_\a
$$
Then we have 
$$
d\phi_\a(t)=(\eta_\a(t)+N(t))\cdot\phi_\a(t).
$$
Taking time derivative of both sides at $t=0$, we have 
$$
d(h\cdot\phi_\a) =(\dot\eta_\a+\dot N)\cdot\phi_\a+(\eta_\a+N)\cdot (h\cdot\phi_\a)
$$
Since $\eta_\a$ is pure imaginary and $N$ is real, we have 
$$
d(h\cdot\ol\phi_\a) =(-\dot\eta_\a+\dot N)\cdot\ol\phi_\a+(-\eta_\a+N)\cdot (h\cdot\ol\phi_\a)
$$
Thus we obtain 
$$
\lan e\cdot\phi_\a, \,\, (d+\eta_\a-N)\cdot (h\cdot\ol\phi_\a)\ran_s =-\lan e\cdot\phi_\a, \,\,(\dot\eta_\a-\dot N)\cdot\ol\phi_\a\ran_s
$$
From Lemma \ref{lanecdotphia, dot Ncdotolphiaran=0}, we have 
$\lan e\cdot\phi_\a, \,\,\dot N\cdot\ol\phi_\a\ran_s=0.$
Hence we obtain 
$$
\lan e\cdot\phi_\a, \,\, (d+\eta_\a-N)\cdot (h\cdot\ol\phi_\a)\ran_s =-\lan e\cdot\phi_\a, \,\,\dot\eta_\a\cdot\ol\phi_\a\ran_s
$$
We decompose $e$ as $e^{1,0}+e^{0,1},$ where $e^{1,0}\in \L_\J$ and $e^{0,1}\in \ol\L_\J.$
We also decompose $\dot\eta_\a=\dot\eta_\a^{1,0}+\dot\eta_\a^{0,1},$ where $\dot\eta_\a^{1,0}\in \L_\J$ and $\dot\eta_\a^{0,1}\in \ol\L_\J.$
Then we have 
\bgn{align*}
-\Im \(i^{-n}\lan e\cdot\phi_\a, \,\, \dot\eta_\a\cdot\ol\phi_\a\ran_s \)=&
-\Im\(i^{-n}\lan e^{0,1}\cdot\phi_\a, \,\, \dot\eta_\a^{1,0}\cdot\ol\phi_\a\ran_s \)\\
=&\Im\( i^{-n}\lan \dot\eta_\a^{0,1}\cdot e^{1,0}\cdot\phi_\a, \,\,\ol\phi_\a\ran_s\)\\
=&\Im\(i^{-n}2\lan e^{0,1}, \,\,\dot\eta_\a^{1,0}\ran_{\tt}\lan \phi_\a, \,\,\ol\phi_\a\ran_s\)
\end{align*}
Since $e$ is real and $\eta_\a$ is pure imaginary, we have 
\bgn{align*}
\lan\dot\eta_\a, \,\,e\ran_{\tt}=&\lan \dot\eta_\a^{1,0}, \,\, e^{0,1}\ran_{\tt}+\lan \dot\eta_\a^{0,1}, \,\,
e^{1,0}\ran_{\tt}\\
=&\lan \dot\eta_\a^{1,0}, \,\, e^{0,1}\ran_{\tt}+\lan -\ol{\dot\eta_\a^{1,0}}, \,\,\ol{e^{0,1}}\ran_{\tt}\\
=&\lan \dot\eta_\a^{1,0}, \,\, e^{0,1}\ran_{\tt}-
\ol{\lan \dot\eta_\a^{1,0}, \,\, e^{0,1}\ran}_{\tt}\\
=&2\sqrt{-1}\,\Im\,\lan \dot\eta_\a^{1,0}, \,\, e^{1,0}\ran_{\tt}
\end{align*}
Since $i^{-n}\lan\phi_\a,\,\,\ol\phi_\a\ran_s =\vol_M,$
we have 
$$
-\Im \(i^{-n}\lan e\cdot\phi_\a, \,\, \dot\eta_\a\cdot\ol\phi_\a\ran_s \)=\Im\,\lan\dot\eta_\a, \,\,e\ran_{\tt}\vol_M
$$
Since $i^{-n}\lan\psi_\a, \,\,\ol\psi_\a\ran_s =\vol_M,$ we have
$$
\Im\, \lan\dot\eta_\a, \,\,e\ran_{\tt}\,\vol_M=\Im\,\lan\dot\eta_\a, \,\,e\ran_{\tt}i^{-n}\lan \psi_\a\, \,\,\ol\psi_\a\ran_s
$$
Then as before, we obtain 
\bgn{align*}
\Im\, \lan\dot\eta_\a, \,\,e\ran_{\tt}\,\vol_M=
&\Im\,\lan\dot\eta_\a, \,\,e\ran_{\tt}i^{-n}\lan \psi_\a\, \,\,\ol\psi_\a\ran_s\\
=&-\Im\, \(i^{-n}\lan e\cdot\psi_\a, \,\,\dot\eta_\a\cdot\ol\psi_\a\ran_s\)
\end{align*}
Since $e$ is a generalized Hamiltonian element,
we have $$e\cdot\psi_\a=-\sqrt{-1}df\cdot\psi_\a.$$ Applying $d\psi_\a=\zeta_\a\cdot\psi_\a,$ we have 
\bgn{align*}
e\cdot\psi_\a=&(-\sqrt{-1}df)\cdot\psi_\a=-\sqrt{-1}\(d(f\psi_\a)-f\zeta_\a\cdot\psi_\a\)\\
=&-\sqrt{-1}(d-\zeta_\a)(f\psi_\a)
\end{align*}
Then we obtain 
\bgn{align}
\Ome_\B(\L_e\J, \,\, \dot\J_h)=&-\Im\(i^{-n}\int_M \lan e\cdot\psi_\a, \,\,\dot\eta_\a\cdot\ol\psi_\a\ran_s\)\\
=&\Im\,\(\int_M   i^{-n+1}\lan (d-\zeta_\a)f\psi_\a,  \,\,\dot\eta_\a\cdot\ol\psi_\a\ran_s \)
\end{align}
%%%%%%%%%%%%%%%%%%%%%%%
We already see that $\phi_\a=e^{ip_{\a,\b}}\phi_\b$. Thus it follows 
$\eta_\a=\eta_\b+id p_{\a,\b}$ where $p_{\a,\b}$ denotes a real function which does not change under small deformations. Thus $\dot\eta_\a=\dot\eta_\b$. (Note that $\eta_\a$ is regarded as a generalized connection form and then infinitesimal deformations of connections is given by 
$\dot\eta_\a$ which is a globally defined section of $\sqrt{-1}(\TT).)$
From the normalization $i^{-n}\lan\psi_\a, \,\,\ol\psi_\a\ran_s=\vol_M$, 
it follows that $f\psi_\a\w\sig(\dot\eta_\a\cdot\ol\psi_\a)=f\psi_\b\w\sig(\dot\eta_\a\cdot\ol\psi_\b)$ for all $\a, \b$.
Thus the component of $(2n-1)$-form of $f\psi_\a\w\sig(\dot\eta_\a\cdot\ol\psi_\a)$ is globally defined.
%%%%%%%%%%%%%%%%%%%%%%%
 Applying the Stokes theorem again, we obtain
 \bgn{align}
 \Ome_\B(\L_e\J, \,\, \dot\J_h)=
 &\Im\(i^{-n+1}\int_M\lan f\psi_\a, \,\, (d+\zeta_\a)\cdot(\dot\eta_\a\cdot\ol\psi_\a)\ran_s\)
 \end{align}
 Then we shall show that 
 $\lan \psi_\a, \,\, (d+\zeta_\a)\cdot(\eta_\a\cdot\ol\psi_\a)\ran_s$
 defines a globally defined $2n$-form on $M.$
 At first,we see 
 \bgn{align*}
 (d+\zeta_\a)\cdot(\eta_\a\cdot\ol\psi_\a)=&d(\eta_\a\cdot\ol\psi_\a)+\eta_\a d\ol\psi_\a
 -\eta_\a d\ol\psi_\a+\zeta_\a\cdot\eta_\a\cdot\ol\psi_\a\\
 =&\L_{\eta_\a}\ol\psi_\a+\eta_\a\cdot\zeta_\a\cdot\ol\psi_\a+\zeta_\a\cdot\eta_\a\cdot\ol\psi_\a\\
 =&\L_{\eta_\a}\ol\psi_\a+2\lan \zeta_\a, \,\eta_\a\ran_{\tt}\ol\psi_\a
\end{align*}
Note that $d\ol\psi_\a=-\zeta_\a\cdot\ol\psi_\a$, since $\zeta_\a$ is pure imaginary.
 Both $\phi_\a$ and $\psi_\a$ satisfy 
 $$
 i^{-n}\lan\phi_\a, \ol\phi_\a\ran_s =i^{-n}\lan \psi_\a, \,\,\ol\psi_\a\ran_s =\vol_M,
 $$
 for all $\a$.
 Thus if $U_\a\cap U_\b$ is not empty,
 we have $\phi_\a=e^{ip_{\a,\b}}\phi_\b$ and $\psi_\a=e^{iq_{\a,\b}}\psi_\b$ for real functions 
 $p_{\a,\b}$ and $q_{\a,\b}.$
 Then we have 
 $\eta_\a=\eta_\b+idp_{\a,\b}$ and $\zeta_\a=\zeta_\b+idq_{\a,\b}.$
 Since $\L_{dp_{\a,\b}}=0$ and $\lan dp_{\a,\b}, \,\, dq_{\a,\b}\ran_{\tt}=0$,  we see
 \bgn{align*}
 \L_{\eta_\a}\ol\psi_\a+2\lan \zeta_\a, \,,\eta_\a\ran_{\tt}\ol\psi_\a  =&\L_{\eta_\b}\ol\psi_\a
 +2\lan \zeta_\b, \,\,\eta_\b\ran_{\tt} \ol\psi_\a+2\lan \zeta_\b, idp_{\a, \b}\ran_{\tt} \ol\psi_\a\\
 +&2\lan idq_{\a,\b},\,\,\eta_\b\ran_{\tt}\ol\psi_\a\\
 =&e^{-iq_{\a,\b}}\(\L_{\eta_\b}\ol\psi_\b+2\lan \zeta_\b, \,\,\eta_\b\ran_{\tt} \ol\psi_\b\)\\
 +&e^{-iq_{\a,\b}}\(2\lan \zeta_\b, idp_{\a, \b}\ran_{\tt} \ol\psi_\b\)
 \end{align*}
 Since $\J_\psi$ is given by a set of locally closed nondegenerate, pure spinor, 
 every $\zeta_\a$ is a one form.
Since $\zeta_\b$ is a one form for all $\b,$
we have $\lan \zeta_\b, idp_{\a, \b}\ran_{\tt}=0.$
Thus it follows
$$
\L_{\eta_\a}\ol\psi_\a+2\lan \zeta_\a, \,,\eta_\a\ran_{\tt}\ol\psi_\a=
e^{-iq_{\a,\b}}\(\L_{\eta_\b}\ol\psi_\b+2\lan \zeta_\b, \,\,\eta_\b\ran_{\tt} \ol\psi_\b\)
$$
Hence we obtain 
$$
\lan \psi_\a, \,\, (d+\zeta_\a)\cdot(\eta_\a\cdot\ol\psi_\a)\ran_s=
\lan \psi_\b, \,\, (d+\zeta_\b)\cdot(\eta_\b\cdot\ol\psi_\b)\ran_s
$$
Hence $\lan \psi_\a, \,\, (d+\zeta_\a)\cdot(\eta_\a\cdot\ol\psi_\a)\ran_s$ defines 
a globally defined $2n$-form on $M$ which is denoted by $i^n\mu(\J)\vol_M,$
\bgn{equation}\label{muJvolM}
i^n\mu(\J)\vol_M:=\lan \psi_\a, \,\, (d+\zeta_\a)\cdot(\eta_\a\cdot\ol\psi_\a)\ran_s,
\end{equation}
where $\mu(\J)$ is a function on $M.$
For an infinitesimal deformation $\{\J_h\}$ in $\B_{\J_\psi}(M),$ 
we have 
\bgn{align}
\Ome_\B(\L_e\J, \,\, \dot\J) =&\Im\(i^{-n+1}\int_M\lan f\psi_\a, \,\, (d+\zeta_\a)\cdot(\dot\eta_\a\cdot\ol\psi_\a)\ran_s\)\\
=&\Im\(i^{-n+1}\int_Mf\frac{d}{dt}\lan\psi_\a, \,\, (d+\zeta_\a)(\eta_\a(t)\cdot\ol\psi_\a)\ran_s\)\Big|_{t=0}\\
=&\Im \(i \int_M f\frac{d}{dt}\mu(\J_t)\vol_M\)\Big|_{t=0}\\
=&\int_M f\frac{d}{dt}\mu(\J_t)\vol_M\Big|_{t=0}
\end{align}
Thus we have 
$\Ome_\B(\L_e\J, \,\, \dot\J)=d\lan\mu, f\ran(\dot\J),$
by using the coupling in terms of the integration over $M.$
Hence 
$\mu: \B_{\J_\psi}(M)\to C^\infty_0(M)$ is a moment map of 
$\B_{\J_\psi}(M)$ for the action of Ham$_{J_\psi}$, which is explicitly given by 
(\ref{muJvolM}) in terms of pure spinors. 
\end{proof}
\bgn{theorem}
Let $X=(M,J, \ome)$ be a compact K\"ahler manifold with a holomorphic Poisson structure $\b$  and
 $(\J_{\b}, \J_\psi)$ a generalized K\"ahler manifold which is given by Poisson deformation.
Then there is a moment map $\mu: \B_{\J_\psi}(M)\to C^\infty_0(M,\R)$ for the action of the group
$\Ham_{\J_\psi}$.
\end{theorem}
\bgn{proof}Since $(\J_{\b}, \J_\psi)$ is a generalized K\"ahler manifold of symplectic type, 
$\J_\psi$ is given by a $d$-closed differential form. 
Thus the result follows form Theorem \ref{existence of moment map}.
\end{proof}
\bgn{definition}
We define the scalar curvature $S(\J)$ of a generalized K\"ahler manifold $(M, \J, \J_\psi)$
to be the moment map $\mu(\J).$
Since the scalar curvature depends on both $\J$ and $\J_\psi$, we also denote the scalar curvature by 
$S(\J, \J_\psi).$
\end{definition}
\section{Lie algebras of generalized complex manifolds and generalized K\"ahler manifolds}
\subsection{The Lie algebra $\frak g_\J$ and the reduced Lie algebra $\frak g_{\red}(\J, \J_\psi)$}
\label{The Lie algebra and the reduced Lie algebra}

Let $\J$ be a generalized complex structure on $M$ which gives the decomposition 
$(\TT)^\C=\L_\J\oplus\ol\L_\J.$
Then we have the Lie algebroid complex:
$$
0\to \w^0\ol \L_\J\overset{\ol\pa_\J}\longrightarrow \w^1\ol\L_\J\overset{\ol\pa_\J}\longrightarrow\w^2\ol\L_\J\overset{\ol\pa_\J}\longrightarrow\cdots\overset{\ol\pa_\J}\longrightarrow\w^n\ol\L_\J\longrightarrow 0
$$
We denote by $H^\bullet(\w^\bullet\ol\L_\J)$ the cohomology groups of the Lie algebroid complex
$\w^\bullet\ol\L_\J.$
\bgn{lemma}\label{the first cohomology}
The first cohomology group $H^1(\w^\bullet\ol\L_\J)$ inherits a Lie algebra structure which is induced from 
the Courant bracket $[\,\,\,,\,\,\,]_{\cou}$
\end{lemma}
\bgn{proof}
Since $\ol\L_\J$ is isotropic, the Courant bracket gives the Lie algebra structure on $\ol\L_\J$, that is, the Jacobi identity holds.
Since $\ol\L_\J$ is a Lie bialgebroid, we have 
$$\ol\pa_\J [e_1, e_2]_{\cou}=[\ol\pa_\J e_1, \,e_2]_{\Sch}+[e_1, \,\,\ol\pa_\J e_2]_{\Sch}$$ for 
$e_1, e_2\in \ol\L_\J,$ where $[\,,\,]_{\Sch}$ denotes the Schouten bracket.
Note that the Courant bracket restricted to $\ol\L_\J$ coincides with 
the Schouten bracket on $\ol\L_\J$.
We also has $\ol\pa_\J[e, f]_{\Sch}=[\ol\pa_\J e, f]_{\Sch}+[e, \ol\pa_\J f]_{\Sch},$
for a function $f$ and $e\in \ol\L_\J.
$
Thus the Courant bracket induces the Lie algebra structure on the first cohomology group 
$H^1(\w^\bullet\ol\L_\J).$
\end{proof}
Then we have 
\bgn{definition}
The Lie algebra $H^1(\w^\bullet\ol\L_\J)$ is denoted by $\frak g_\J$ which is called the Lie algebra of automorphisms of $(M, \J)$.
\end{definition}
We have the following lemma:
\bgn{lemma}\label{e1,e2=dlane1,e2ran} 
If $e_1, e_2\in \ol\L_\J$ and  satisfy $\ol\pa_\J e_1=\ol\pa_\J e_2=0.$
Then the Courant bracket of $e_1$ and the conjugate $\ol e_2$ is given by 
$$
[e_1, \,\,\ol e_2]_{\cou}=(\pa_\J-\ol\pa_\J)\lan e_1, \,\ol e_2\ran_{\tt}
$$
\end{lemma}
\bgn{proof}
Let $\{A, B\}$ be the anti-bracket $AB+BA$ for operators $A,B$ in general.
Then we have the super-Jacobi identity
$$
[A, \,\{B,C\}]=[\{A, B\}, C]+[\{A, C\}, B].
$$
Thus we have 
$$
[d, \{ e_1, \ol e_2\}]=[\{d, e_1\}, \ol e_2]+[\{d, \ol e_2\}, e_1]
$$
Since $\{e_1, \ol e_2\}=2\lan e_1, \ol e_2\ran_{\tt},$
we obtain
$$
2d\lan e_1, \ol e_2\ran_{\tt}=[\{d, e_1\}, \ol e_2]+[\{d, \ol e_2\}, e_1]
$$
From the definition of the Courant bracket, 
we have 
$[e_1, \ol e_2]_{\cou}=\frac12 [\{d, e_1\}, \ol e_2]-\frac12[\{d, \ol e_2\}, e_1].$
Since $\ol\pa_\J e_1=\ol\pa_\J e_2=0,$
we also have 
$$
[\{d, e_1\}, \ol e_2]+[\{d, \ol e_2\}, e_1]=[\{\pa_\J, e_1\}, \ol e_2]+[\{\ol\pa_\J, \ol e_2\}, e_1]
$$
Since $\pa_\J \ol e_2=\{\pa_\J,\ol e_2\}=0$ and $\ol\pa_\J e_1=\{\ol\pa_\J ,e_1\}=0, $
applying the super-Jacobi identity again, 
we have  
\bgn{align}
[\{\pa_\J, e_1\}, \ol e_2]-[\{\ol\pa_\J, \ol e_2\}, e_1]=
&\pa_\J\{e_1, \ol e_2\}-\ol\pa_\J \{\ol e_2, e_1\}\\
=&2(\pa_\J-\ol\pa_\J)\lan e_1, \ol e_2\ran_{\tt}.
\end{align}
Hence we have 
$$
[e_1, \,\ol e_2]_{\cou }=(\pa_\J-\ol\pa_\J)\lan e_1, \ol e_2\ran_{\tt}
$$
\end{proof}
We define a map $F: \ol\L_\J \to (\TT)$ by taking the real part of $e\in \ol\L_\J$, 
\bgn{equation}\label{f(e):=e+ol e}
F(e):=e+ \ol e
\end{equation}
Then $F$ restricted to $\ker\ol\pa_\J$ yields the map from $\ker\ol\pa_\J$ to the real part
$(\ker\ol\pa_\J +\ker\pa_\J)^\R.$
Taking the quotient, we have the map from 
$\frak g_\J$ to $(\ker\ol\pa_\J +\ker\pa_\J)^\R/(\Im\,\ol\pa_\J+\Im\,\pa_\J)^\R.$
By the abuse of notation, we denote by $F
$ the map to the quotient.
Then we have 
\bgn{proposition}
The quotient $(\ker\ol\pa_\J +\ker\pa_\J)^\R/(\Im\,\ol\pa_\J+\Im\,\pa_\J)^\R$
is a Lie algebra with respect to the Courant bracket and
$$
F: \frak g_\J \to (\ker\ol\pa_\J +\ker\pa_\J)^\R/(\Im\,\ol\pa_\J+\Im\,\pa_\J)^\R
$$
is an isomorphism between Lie algebras.
\end{proposition}
\bgn{proof}
For $e_1, e_2\in \ker \ol\pa_\J,$  from Lemma \ref{e1,e2=dlane1,e2ran} and taking the complex conjugate, one has 
\bgn{align}
[e_1+\ol e_1, \, e_2+\ol e_2]_{\cou}=&[e_1, e_2]_{\cou}+[\ol e_1, \ol e_2]_{\cou}+[e_1, \ol e_2]_{\cou}+[\ol e_1, e_2]_{\cou}\\
=&[e_1, e_2]_{\cou}+[\ol e_1, \ol e_2]_{\cou}\\
+&(\pa_\J-\ol\pa_\J)\(\lan e_1,\,\,\ol e_2\ran_{\tt}-\lan \ol e_1, \,\,e_2\ran_{\tt}\),
\end{align}
since $(\pa_\J-\ol\pa_\J)\(\lan e_1,\,\,\ol e_2\ran_{\tt}-\lan \ol e_1, \,\,e_2\ran_{\tt}\)
\in 
(\Im\,\ol\pa_\J+\Im\,\pa_\J)^\R.$
Hence $(\ker\ol\pa_\J +\ker\pa_\J)^\R/(\Im\,\ol\pa_\J+\Im\,\pa_\J)^\R$
is a Lie algebra with respect to the Courant bracket.
Since we see 
$$[F(e_1), \,\, F(e_2)]_{\cou}\equiv F([e_1, e_2]_{\cou}) \quad \text{   \rm mod } 
(\Im\,\ol\pa_\J+\Im\,\pa_\J)^\R,$$
thus $F$ is an isomorphism between Lie algebras.
\end{proof}
Let $(M, \J, \J_\psi)$ be a generalized K\"ahler manifold. 
Then we define a subspace $\wtil{\frak g}_{\red}(\J, \J_\psi)$ of $\ker\ol\pa_\J$ by 
$$
\wtil{\frak g}_{\red}(\J, \J_\psi):=
\{\, \J_\psi (\ol\pa_{\ss\J} u)\,|\, \ol\pa_{\ss\J}\J_\psi \ol\pa_{\ss\J}u=0,\,\, u\in C^\infty(M, \C)\, \}
\subset\ol\L_\J
$$
For simplicity, we also denote by $\wtil{\frak g}_{\red}$ the subspace $\wtil{\frak g}_{\red}(\J, \J_\psi)$ .
Since $\wtil{\frak g}_{\red}$ is a subspace of $\ker\ol\pa_{\ss\J}$, we have the following diagram:
$$
\xymatrix{
&\wtil{\frak g}_{\red}\ar@{->}[dr]^j\ar@{->}[r]^i&\ker \ol\pa_{\ss\J}\ar@{->}[d]\\
&&\frak g_\J:=\ker\ol\pa_{\ss\J}/\Im \ol\pa_{\ss\J}
}
$$
Then we have 
\bgn{proposition}
If $M$ is compact, then the map $j: \wtil{\frak g}_{\red}\to\frak g_\J$ is injective.
\end{proposition}
\bgn{proof}
It suffices to show that the intersection $\wtil{\frak g}_{\red}\cap \Im\,\ol\pa_{\ss\J}=\{0\}.$
We assume that there exist two functions $u$ and $v$ such that 
$$
\J_\psi (\ol\pa_{\ss\J} u)=\ol\pa_{\ss\J}v.
$$
Since $\ol\pa_{\ss\J}u-\sqrt{-1}\J_\psi(\ol\pa_{\ss\J}u)\in \ol\L_{\J_\psi},$
we have $\ol\pa_{\ss\J}u-\sqrt{-1}\,\ol\pa_{\ss\J}v=\ol\pa_{\ss\J}(u-\sqrt{-1}\,v)\in \ol\L_{\J_\psi}.$
We have the decomposition $\ol\L_\J=\ol\L_\J^+\oplus\ol\L_\J^-$, where $\ol\L_\J^+=\ol\L_\J\cap \ol\L_{\J_\psi}$ and $\ol\L_\J^-=\ol\L_\J\cap \L_{\J_\psi}.$
Thus $\ol\pa_\J (u-\sqrt{-1}v)\in \ol\L_\J^+.$
Hence 
\bgn{align}\label{olpaJ-u-sqrt}\ol\pa_\J^-(u-\sqrt{-1}v)&=0,\end{align}
where $\ol\pa_\J=\ol\pa_\J^++\ol\pa_\J^-.$
By $\ol\pa_{\ss\J}u+\sqrt{-1}\J_\psi(\ol\pa_{\ss\J}u)\in \L_{\J_\psi},$
we have $\ol\pa_{\ss\J}u+\sqrt{-1}\,\ol\pa_{\ss\J}v=\ol\pa_{\ss\J}(u+\sqrt{-1}\,v)\in \L_{\J_\psi}.$
Thus we also have $\ol\pa_\J (u+\sqrt{-1}v)\in \ol\L_\J^-.$
Hence 
\bgn{equation}\label{olpaJJ+squrt}\ol\pa_\J^+(u+\sqrt{-1}v)=0.
\end{equation}
Thus we have 
\bgn{align*}
&\ol\pa_\J^-(u-\sqrt{-1}v)\cdot\ol\psi =\ol\del_-\((u-\sqrt{-1}v)\ol\psi\)=0\\
&\ol\pa_\J^+(u+\sqrt{-1}v)\cdot\psi=\ol\del_+\((u+\sqrt{-1}v)\psi\)=0
\end{align*}
Since the generalized K\"ahler identity, the Laplacian $(\ol\del_\pm)^*\ol\del_\pm +\ol\del_\pm(\ol\del_\pm)^*$ of the operator $\ol\del_\pm$ is $\frac 14 \trian$, where $\trian$ denotes the ordinary Laplacian $dd^*+d^*d.$
Since $(u\pm\sqrt{-1}v)$ is a function,
then it follows from (\ref{olpaJ-u-sqrt}) and (\ref{olpaJJ+squrt})that $\trian (u\pm\sqrt{-1}v)=0$.
Thus $u+\sqrt{-1}v$ and $u-\sqrt{-1}v$ are constants. 
Thus $u$ and $v$ are constants also.
Hence  
$\J_\psi (\ol\pa_{\ss\J} u)=\ol\pa_{\ss\J}v=0.$
Thus we have $\frak g_{\red}\cap \Im\,\ol\pa_{\ss\J}=\{0\}.$
\end{proof}
\bgn{definition}
We define $\frak g_{\red}$ to be the image $j(\wtil{\frak g}_{\red})$ in $\frak g_\J.$
\end{definition}
\bgn{proposition}\label{wtilfrakgredsubsetH1}
${\frak g}_{\red}\subset \frak g_\J:=H^1(\w^\bullet\ol\L_\J)$
is a Lie subalgebra.
\end{proposition}
\bgn{proof}
Since
$\J_\psi$ is integrable, the Nijenhuis tensor vanishes,
\bgn{align}\label{[Jpsi e1, Jpsi e2]}
[\J_\psi e_1, \J_\psi e_2]_{\cou}=[e_1, e_2]_{\cou}+\J_\psi[\J_\psi e_1, \, e_2]_{\cou}+\J_\psi[e_1, \J_\psi e_2]_{\cou},
\end{align}
where $e_1, e_2\in \ol\L_\J.$
For simplicity, we denote by $\ol\pa$ the operator $\ol\pa_{\J}.$
For $u, v\in C^\infty(M, \C)$, we assume 
$\J_\psi(\ol\pa u)$ and $ \J_\psi(\ol\pa v)\in \ol\L_\J$ satisfy
$\ol\pa\J_\psi(\ol\pa u)=0$ and $ \ol\pa\J_\psi(\ol\pa v)=0,$ respectively.
Then from (\ref{[Jpsi e1, Jpsi e2]}), we have 
\bgn{align}
[\J_\psi(\ol\pa u), \,\J_\psi(\ol\pa v)]_{\cou}=&[\ol\pa u,\,\ol\pa v]_{\cou}+
\J_\psi[\J_\psi(\ol\pa u),\,\ol\pa v]_{\cou}+\J_\psi [\ol\pa u,\,\J_\psi(\ol\pa v)]_{\cou}
\end{align}
Since $\ol\L_\J$ is a Lie bialgebroid,
we have 
$ \ol\pa [e_1, e_2]_{\cou}=[\ol\pa e_1, e_2]_{\cou}+[e_1, \ol\pa e_2]_{\cou}$, for $e_1, e_2\in \ol\L_\J.$
Thus we have 
$[\ol\pa u,\,\ol\pa v]_{\cou}=\ol\pa[u, \ol\pa v]_{\cou}$ which vanishes as an element of $H^1(\w^\bullet\ol\L_\J)$. 
From our assumption $\ol\pa\J_\psi(\ol\pa u)=0, \,\,\ol\pa\J_\psi(\ol\pa v)=0$, we have 
\bgn{align}
\J_\psi[\J_\psi(\ol\pa u),\,\ol\pa v]_{\cou}=&\J_\psi\ol\pa [\J_\psi(\ol\pa u), \,\, v]_{\cou}\\
\J_\psi [\ol\pa u,\,\J_\psi(\ol\pa v)]_{\cou}=&\J_\psi\ol\pa[u, \,\,\J_\psi(\ol\pa v)]_{\cou}
\end{align}
Thus we obtain
$$
[\J_\psi(\ol\pa u), \,\J_\psi(\ol\pa v)]_{\cou}=\ol\pa[u, \ol\pa v]_{\cou}
+\J_\psi\ol\pa\{ [\J_\psi(\ol\pa u), \,\, v]_{\cou}+[u, \,\,\J_\psi(\ol\pa v)]_{\cou}
\},
$$
where $ [\J_\psi(\ol\pa u), \,\, v]_{\cou}:=\L_{\J_\psi(\ol\pa u)}v$ and 
$[u, \,\,\J_\psi(\ol\pa v)]_{\cou}:=-\L_{\J_\psi(\ol\pa v)}u$ and \\
$ [\J_\psi(\ol\pa u), \,\, v]_{\cou}+[u, \,\,\J_\psi(\ol\pa v)]_{\cou}\in C^\infty(M, \C)$.  
We denote by $\{u, v\}_{\J,\psi}$ the complex function 
$ [\J_\psi(\ol\pa u), \,\, v]_{\cou}+[u, \,\,\J_\psi(\ol\pa v)]_{\cou}.$
Then we have 
$$
[\J_\psi(\ol\pa u), \,\J_\psi(\ol\pa v)]_{\cou}\cong \J_\psi(\ol\pa\{u, v\}_{\J,\psi})
\in H^1(\w^\bullet\ol\L_\psi)
$$
Hence the result follows.
\end{proof}
Then
$\frak g_{\red}$ is called the Lie algebra of reduced automorphisms of $(M, \J, \J_\psi).$
\subsection{The real Lie algebra $\frak g_{\red}^{\Re}$ of the Lie algebra of the reduced automorphisms}
\label{the rela Lie algebra frak gredR}
In this section, we assume that $(M, \J, \J_\psi)$ is a compact generalized K\"ahler manifold of symplectic type, i.e., $\J_\psi$ is given by a $d$-closed, nondegenerate, pure spinor
$\psi=e^{B-\sqrt{-\ome}}$, where $\ome$ denotes a real symplectic $2$-form on $M$ and $B$ is a real $d$-closed $2$-form on $M.$
Since the map $j$ is injective, $j$ gives an isomorphism between $\frak g_{\red}$ and $\wtil{\frak g}_{\red}.$
Thus  $\frak g_{\red}$ is identified with $ \wtil{\frak g}_{\red}.$
\bgn{definition}
Consider a real Lie subalgebra $\frak g_{\red}^{\Re}$ of $\frak g_{\red}$ by 
$$
\frak g_{\red}^{\Re}:=\{\, \J_\psi\ol\pa_{\ss\J} u\in \frak g_{\red}\, |\, u\in C^\infty(M, \R)\, \}
$$
\end{definition}
A Lie algebra is called {\it a reductive Lie algebra} if the radical of the Lie algebra equals the center, where the radical is the maximal solvable ideal.
A reductive Lie algebra is the direct sum of 
a semisimple Lie algebra and an abelian Lie algebra. 
It is known that a Lie algebra is reductive if the associated Lie group of the Lie algebra is a compact Lie group. 
\bgn{proposition}\label{compact Lie algebra}
The real sub Lie algebra $\frak g_{\red}^{\Re}$ is a reductive Lie algebra if a compact generalized K\"ahler manifold
$(M, \J, \J_\psi)$ is of symplectic type.
\end{proposition}
\bgn{proof}
By using the $B$-field transformation, Proposition \ref{compact Lie algebra} reduces to
 the case $B=0$. 
 Thus it suffices to show Proposition in the case $B=0.$ 
By using the map $F$ as in (\ref{f(e):=e+ol e}), 
it follows from $u\in C^\infty(M,\R)$ that  $F(\J_\psi \ol\pa_{\ss\J}u) =\J_\psi du.$
Since $B=0,$ $\J_\psi du$ is an ordinary Hamiltonian vector field with respect to $\ome.$
Thus the real Lie algebra 
$\frak g_{\red}^{\Re}$ is isomorphic to $F(\frak g_{\red}^{\Re})$ which is a subgroup of the Lie algebra of Hamiltonian vector fields
$$
\{\, \J_\psi du \in T_M\, | u\in C^\infty(M,\R)\,\}
$$
Since a Hamiltonian vector field acts on $M$ preserving $\J_\psi$ and 
$\frak g_{\red}^{\Re}$ also preserves $\J$, 
it follows that $\frak g_{\red}^{\Re}$ preserves the generalized metric $G$ of 
$(M, \J, \J_\psi).$
The generalized metric $G$ consists of a Riemannian metric $g$ and a $2$-form $b$ which satisfies 
$d^c_{I_+}\ome_{I_+}=-d^c_{I_-}\ome_{I_-}=db.$
Thus a Hamiltonian vector field $\J_\psi du \in F(\frak g_{\red}^{\Re})$
is a Killing vector field with respect to $g$ which preserves $b.$ 
Since
$\L_e\J=0$ is equivalent to
 $(\ol\pa_{\ss\J}e)=0$ for $e\in \ol\L_{\ss\J}$, if a Hamiltonian vector field $\J_{\psi}du$ is a Killing vector field preserving $b$, then $\J_\psi du\in \frak g_{\red}^{\Re}.$
Thus $\J_\psi du\in \frak g_{\red}^{\Re}$ if and only if 
$\J_\psi du\in \frak g_{\red}^{\Re}$ is a Killing vector field which preserves $b.$
We denote by $G_{\red}^{\Re}$ the associated Lie group with $\frak g_{\red}^{\Re}$. 
Then $G_{\red}^{\Re}$ is a subgroup of the isometry group Isom$(M, g)$ of the Riemannian manifold $(M, g)$.
It is know that Isom$(M,g)$ is a compact Lie group of finite dimension.
Let Symp$_0(M,\ome)$ be the identity component of diffeomorphisms which preserves $\ome.$
We denote by Ham$(M,\ome)$ the group of Hamiltonian diffeomorphisms.
Then the following theorem is known as the Flux conjecture which is affirmatively solved.
\bgn{theorem}\label{Ono}\text{\rm \cite{OK_2006}}
$\text{\rm Ham}(M,\ome)$ is a closed subgroup of Symp$_0(M,\ome)$ with respect to $C^1$-topology.
\end{theorem}
Thus it follow from Theorem \ref{Ono} that 
the intersection Ham$(M,\ome)\,\cap $Isom$(M,g)$ is a compact Lie group.
The group $G_{\red}^{\Re}$ is a subgroup of Ham$(M,\ome)\cap\,$Isom$(M,g)$ which preserves $b.$ 
Let $\{f_i\}$ be a set of $G_{\red}^{\Re}$ which converges to a function $f_\infty\in$ Ham$(M,\ome)\cap $Isom$(M,g)$ with respect to $C^r$-topology
for $r\geq 1.$
Then we have $\lim_{i\to\infty}f_i^*b=f_\infty^*b$. 
Since $f_i^*b=b$, we have $f_\infty^*b=b$. 
Thus $G_{\red}^{\Re}$ is a closed  subgroup of a compact Lie group Ham$(M,\ome)\cap $Isom$(M,g)$, which is also
a compact Lie group.
Hence $\frak g_{\red}^{\Re}$ is a Lie algebra of a compact Lie group. 
Thus $\frak g_{\red}^{\Re}$ is a reductive Lie algebra.
\end{proof}
\subsection{Reductivety of $\frak g_{\red}$}\label{Reductiveity of frak gred}
Let $\frak g_{\red}$ be the Lie algebra of reduced automorphisms of $(M, \J, \J_\psi)$ as before.
A complex function $u$ is denotes by $u_{\Re}+\sqrt{-1}u_{\Im}$, where 
$u_{\Re}$ is the real part of $u$ and $u_{\Im}$ is the imaginary part of $u.$
We consider the following condition (\ref{uRe uIm}) on a generalized K\"ahler manifold $(M,\J,\J_\psi)$ :\\
\\ \noindent
If a complex function $u$ satisfies $\ol\pa \J_\psi \ol\pa u =0, $
then both $u_{\Re}$ and $u_{\Im}$ also satisfy 
\bgn{equation}\label{uRe uIm}
\qquad \qquad \ol\pa\J_\psi \ol\pa u_{\Re}=0, \quad \ol\pa\J_\psi\ol\pa u_{\Im}=0
\end{equation}
\bgn{proposition}\label{complexification}
Let $(M, \J,\J_\psi)$ be a compact generalized K\"ahler manifold of symplectic type which satisfies the condition
(\ref{uRe uIm}).
Then $\frak g_{\red}$ is the complexification of the Lie algebra of a compact Lie group, that is, 
$\frak g_{\red}$ is a reductive Lie algebra.
\end{proposition}
\bgn{proof}
The condition (\ref{uRe uIm}) implies that $\frak g_{\red}$ is the complexification of $\frak g_{\red}^{\Re}.$
Then the result follows from Proposition \ref{compact Lie algebra}.
\end{proof}
\subsection{The structure theorem of the Lie algebra ${\frak g}_{\ss\red}$ and the Lie algebra ${\frak g}_\J$}
\label{The structure theorem of the Lie algebra}
Let $(M, \J, \J_\psi)$ be a compact generalized K\"ahler  manifold of symplectic type. 
Then we have the decomposition $\ol \L_\J =\ol\L_\J^+\oplus \ol\L_\J^-$
 and $\ol\pa_\J=\ol\pa_++\ol\pa_-$. 
 We define $\ol\L^{p,q}:=\w^p\ol\L_\J^+ \otimes \w^q\ol\L_\J^-.$
Then we have the double complex $(\ol\L^{\bullet, \bullet}, \ol\pa_+, \ol\pa_-)$. 
In generalized K\"ahler manifold, the space of differential forms is decomposed into $\oplus_{p,q}U^{p,q},$
where $-n\leq p+q\leq n$ and $-n\leq -p+q\leq n$.
The exterior derivative $d$ is decomposed into $\del_++\del_-+\ol\del_++\ol\del_-$ and 
it is known that the generalized K\"ahler identity does holds. 
Then the double complex $(U^{\bullet, \bullet},\ol\del_+, \ol\del_-)$ defines the cohomology groups $H^{p,q}(M, \J, \J_\psi)$ and 
the generalized Hodge decomposition holds: 
$$
\oplus_{i=0}^{2n}H^i(M, \C)=\oplus_{p,q}H^{p,q}(M, \J, \J_\psi)
$$
The isotropic space $\ol\L_\J^{p,q}$ acts on $\psi$ by the Spin action which is given by the interior product and the exterior product. 
Then we see  $U^{p, -n+q}=\ol\L_\J^{p,q}\cdot\psi$ and $U^{p, n-q}=\ol\L_\J^{p,q}\cdot\ol\psi,$
where $\ol\psi$ denotes the complex conjugate of $\psi.$
We denote by $H^{\odd}(M,\C)$ the direct sum $\oplus_{i=0}^{n-1}H^{2i+1}(M,\C)$ of the de Rham cohomology groups of odd degree. 
Let $[a]$ be a class in $\frak g_\J=H^1(\w^\bullet\ol\L_\J)$.
Then the representative $a=a_++a_-$ is a $\ol\pa_\J$-closed element of $\ol\L_\J=\ol\L_++\ol\L_-$, where $a_+\in \ol\L_+$ and 
$a_-\in\ol\L_-.$ 
The condition $\ol\pa_\J a=0$ is equivalent to $\ol\pa_+a_+=0, \ol\pa_-a_-=0$ and $\ol\pa_-a_++\ol\pa_+a_-=0.$
Then $a_+$ acts on $\psi$ to obtain $a_+\cdot\psi\in U^{1,-n+1}$ and $a_-$ also acts on $\ol\psi$
to get $a_-\cdot\ol\psi\in U^{1,n-1}.$
Since $\psi$ is $d$-closed, we obtain 
$$
\ol\del_+(a_+\cdot\psi)=(\ol\pa_+a_+)\cdot\psi =0, \qquad \ol\del_-(a_-\cdot\ol\psi)=(\ol\pa_- a_-)\cdot\ol\psi=0.
$$
If $a=\ol\pa_\J u =\ol\pa_+ u +\ol\pa_-u,$ then we have 
$$a_+\cdot\psi =(\ol\pa_+u)\cdot\psi =\ol\del_+(u\psi),\qquad
a_-\cdot\ol\psi =(\ol\pa_- u)\cdot\ol\psi=\ol\del_-(u\ol\psi).$$
Thus $a_+\psi$ defines a class $[a_+\psi]\in H^{1,-n+1}(M,\J,\J_\psi)$ and $a_-$ also defines a class $[a_-\cdot\ol\psi]\in H^{1,n-1}(M,\J,\J_\psi)$. 
Thus we have a map $$ \frak g_\J \to H^{1,-n+1}(M,\J,\J_\psi)\oplus  H^{1,n-1}(M,\J,\J_\psi).$$
Since $\psi$ is a differential form $e^{b-\sqrt{-1}\ome}$, it follows from the generalized Hodge decomposition that 
$H^{1,-n+1}(M,\J,\J_\psi)\oplus  H^{1,n-1}(M,\J,\J_\psi)$ is isomorphic to $H^{1}(M,\C).$
Im fact, $(a_+\cdot\phi)\w e^{-b+\sqrt{-1}\ome}+(a_-\cdot\ol\phi)\w e^{-b-\sqrt{-1}\ome}$ gives a representative of $H^1(M,\C)$.
Thus we obtain a map $j : \frak g_\J \to H^{1}(M, \C).$
\bgn{theorem}[the structure theorem of $\frak g_\J$]\label{the structure theorem of frakgJ}
Then the following exact sequence of Lie algebras holds
$$
0\to \frak g_{\red}\overset i\to \frak g_\J\overset j\to H^1(M,\C)\to 0,
$$
where $\frak g_{\red}$ is the Lie algebra of reduced automorphisms of $(M, \J, \J_\psi)$, which is a Lie subalgebra of $\frak g_\J$, and $\frak g_{\red}$ has 
a Lie subalgebra $\frak g_{\red}^{\Re}$ which is a real reductive Lie algebra.
Moreover, $H^1(M,\C)$ is a commutative Lie algebra.
\end{theorem}
\bgn{proof}
First we shall show that $\ker j =i(\frak g_\J).$
If a class $[a]=[a_++a_-]\in \frak g_\J$ satisfies $j([a])=0\in H^1(M,\C),$ then there exist two functions 
$u, v$ such that
$$
a_+\cdot\psi =\ol\del_+(u\psi)=(\ol\pa_+u)\cdot\psi, \qquad a_-\cdot\ol\psi =\ol\del_- (v\ol\psi)=(\ol\pa_-v)\cdot\ol\psi.
$$
Thus we see $a_+=\ol\pa_+ u, \,\, a_-=\ol\pa_- v$.
Since $\ol\pa_\pm =\frac12(\ol\pa_\J \mp\sqrt{-1}\J_\psi \ol\pa_\J),$ we have 
\bgn{align*}
a=&a_++a_-=\frac12(\ol\pa_\J -\sqrt{-1}\J_\psi \ol\pa_\J)u+ \frac12(\ol\pa_\J +\sqrt{-1}\J_\psi \ol\pa_\J)v\\
=&\frac12\ol\pa_\J (u+v)-\frac12\sqrt{-1}\J_\psi\ol\pa_\J (u-v)
\end{align*}
Thus the class $[a]$ is represented by $-\frac12\sqrt{-1}\J_\psi\ol\pa_\J (u-v)\in \frak g_{\red}.$
Hence $\ker j\subset i(\frak g_{\red}).$
Conversely, a class $i(\frak g_{\red})$ is represented by $\J_\psi \ol\pa u =\sqrt{-1}\,\ol\pa_+ u -\sqrt{\-1}\,\ol\pa_- u.$ Then we see $i(\frak g_{\red})\subset \ker j.$
Hence $\ker j =i(\frak g_{\red}).$
From the generalized Hodge decomposition theorem, it follows that $j$ is surjective. 
By using the generalized $\pa_\J\ol\pa_\J$-lemma, we obtain 
$[a, b]_{co}\in \frak g_0 $ for all $a, b \in \ol\L_\J$ satisfying $\ol\pa_\J a =\ol\pa_\J b=0$.
Thus the quotient $\frak g_\J/\frak g_0 =H^1(M,\C)$ is a commutative  Lie algebra.
\end{proof}
\bgn{corollary}\label{cor hodd=0}
Let $(M, \J, \J_\psi)$ be a compact generalized K\"ahler manifold of symplectic type.
If $H^{1}(M,\C)=0$, then $\frak g_\J \cong \frak g_{\red}.$
\end{corollary}
\bgn{proof}
The result follows from Theorem \ref{the structure theorem of frakgJ}.
\end{proof}

\subsection{The Lie algebras ${\frak g}_\J$ and ${\frak g}_{\red}$ of generalized K\"ahler manifolds which are given by Poisson deformations}
\label{The Lie algebra poisson deformations}
Let $X=(M, J)$ be a compact complex manifold with a K\"ahler form $\ome$ and 
$\b$ a holomorphic Poisson structure.
We assume that $H^1(X, {\cal O})=0$ in this section.
We denote by $\{(M, \J_{\b t}, \J_{\psi_t})\}$ a family of generalized K\"ahler manifolds which is given by 
Poisson deformations, where $\psi_t=e^{b_t-\sqrt{-1}\ome_t}$ is the $d$-closed, nondegenerate, pure spinor, where $t$ is a parameter of deformations. Note that $\ome_t$ is a symplectic form which is not of type $(1,1)$ with respect to the ordinary complex structure $J.$
Then we have the Lie algebroid complex :
$$
0\to C^\infty(M, \C)\overset{\ol\pa_{\J_{\b t}}}\arrow \ol \L_{\J_{\b t}}\overset{\ol\pa_{\J_{\b t}}}\arrow\w^2\ol\L_{\J_{\b t}}
\overset{\ol\pa_{\J_{\b t}}}\arrow\cdots
$$
Then we have the Lie algebra $H^1(\w^\bullet\ol\L_{\J_{\b t}})$ 
and we see that the Lie algebra $H^1(\w^\bullet\ol\L_{\J_{\b t}})$ does not depend on $t\neq 0.
$
For simplicity, we denote by $\frak g_{\J_\b}$ the Lie algebra $H^1(\w^\bullet\ol\L_{\J_{\b t}})$.
\bgn{proposition}\label{Poisson vector fields}
We assume that $H^1(X, {\cal O})=0.$
Then the Lie algebra ${\frak g}_{\J_{\b }}$ of automorphisms with respect to $\J_{\b }$ is given by the Lie algebra of holomorphic vector fields preserving the holomorphic Poisson structure $\b$, i.e.,
$$
{\frak g}_{\J_{\b }}=\{ V\in H^0(X, T^{1,0}_J)\, |\, \L_V\b=0\, \}.
$$
%Further, $\wtil{\frak g}(\J_{\b t})=\wtil{\frak g}(\J_\b)$ for $t\neq 0.$
\end{proposition}
\bgn{proof}
For simplicity, we denote by ${\frak g}_\b$ the Lie algebra ${\frak g}_{\J_{\b }}\cong H^1(\w^\bullet\ol\L_{\J_\b}).$
The cohomology $H^1(\w^\bullet\ol\L_{\J_\b})$ is the total cohomology of the double complex $(\w^pT^{1,0}_J\otimes
\w^{0,q}_\b, \del_\b, \ol\pa_\b)$, 
where $\w^{0,1}_\b=\{\t+[\ol\b, \t]\, |\, \t\in \w^{0,1}\, \}$ and $\ol\pa_\b=e^{-\ol\b}\circ \ol\pa\circ e^{\ol\b}$
and $\w^{0,q}_\b=\w^q(\w^{0,1}_\b)$ and 
$\del_\b : \w^p T^{1,0}_J \to \w^{p+1}T^{1,0}_J$ denotes the Poisson complex. 
The complex $(\w^pT^{1,0}_J\otimes
\w^{0,q}_\b,  \ol\pa_\b)$ for each $p$ is quasi-isomorphic to 
the ordinary Dolbeault complex 
$(\w^pT^{1,0}_J\otimes\w^{0,q},  \ol\pa).$
Thus taking the cohomologies by using $\ol\pa_\b$ at first,
we have the $E_1$-terms in terms of the ordinary Dolbeault cohomology groups, 
$$
E^{p,q}_1=H^{0,q}(X, \w^pT^{1,0}_J)
$$
Since $H^1(X, {\cal O})=0,$ we have 
$$
E_1^{0,1}=0, \qquad E^{1,0}_1=H^0(X, T^{1,0}_J)
$$
Thus the $E_2$-terms are given by 
$$
E_2^{0,1}=0, \qquad 
E_2^{1,0}=\Ker \del_\b: H^0(X, T^{1,0}_J)\to H^0(X, \w^2T^{1,0}_J)
$$
Thus the total cohomology $H^1(\w^\bullet\ol\L_{\J_\b})$
is $E_2^{1,0}.$
Since $\del_\b V$ is given by the Lie derivative $\L_V\b$ of $\b$ by $V$, we have 
$H^1(\w^\bullet\ol\L_{\J_\b})=\{\, V\in H^0(X, T^{1,0}_J)\, \,|\,\, \L_V\b=0\,\, \}.$
Hence we see that ${\frak g}_\b =\{\, V\in H^0(X, T^{1,0}_J)\, \,|\,\, \L_V\b=0\,\, \}$.
Thus the result follows.
\end{proof}

 We also have the Lie algebra of the reduced automorphisms ${\frak g}_{\red}(\J_{\b t}, \J_{\psi, t})
\subset 
H^1(\w^\bullet\ol\L_{\J_\b t})$ of a generalized K\"ahler manifold $(M, \J_{\b t}, \J_{\psi, t})$ as in Section \ref{The Lie algebra and the reduced Lie algebra}. 
\bgn{proposition}\label{reduced auto coincides with}
We assume that $H^{1}(M, \C)=0.$
Then the Lie algebra ${\frak g}_{\J_\b}$ of automorphisms coincides with the Lie algebra $\frak g_{\red}:={\frak g}_{\ss \red}(\J_{\b t}, \J_{\psi, t})$ of the reduced automorphisms 
of a generalized K\"ahler manifold $(M, \J_{\b t}, \J_{\psi_t}).$ 
\end{proposition}
\bgn{proof}
The result follows from Corollary \ref{cor hodd=0}.
\end{proof}
\section{The Lie algebra of automorphisms of generalized K\"ahler manifolds with constant scalar curvature
(Matsushima -Lichnerowicz type theorems)}
\label{The Lie algebra of cscGK}
Let
$(M,\J, \J_\psi)$ be a generalized K\"ahler manifold of symplectic type, where 
$\psi=e^{b-\sqrt{-1}\ome}$ and $\ome$ is a symplectic form.
Then the generalized metric $G=-\J\circ\J_\psi$ gives a Hermitian metric on 
$\w^\bullet\ol\L_\J$. 
The operator $\ol\pa_\J: \w^\bullet\ol\L_J\to\w^{\bullet+1}\ol\L_\J$ together with 
$\J_\psi:\ol\L_\J\to \ol\L_\J$ gives a second order differential operator 
$\ol\pa_\J\J_\psi \ol\pa_\J :C^\infty(M, \C) \to \w^2\ol\L_\J.$
The adjoint operator of $\ol\pa_\J\J_\psi \ol\pa_\J$ is denoted by $(\ol\pa_\J\J_\psi\ol\pa_\J)^*: \w^2\ol\L_\J\to C^\infty(M, \C)$.
Then we define the fourth order differential operator $L:C^\infty(M, \C)\to C^\infty(M, \C)$
by the composition
\bgn{equation}\label{L=(olpaJpsiolpa}
L=(\ol\pa_\J\J_\psi\ol\pa_\J)^*\circ(\ol\pa_\J\J_\psi\ol\pa_\J).
\end{equation}
We denote by $\ol L$ the complex conjugate of the operator $L$, i.e., 
\bgn{equation}\label{ol L=(paJpsi}
\ol L =(\pa_\J\J_\psi\pa_\J)^*\circ(\pa_\J\J_\psi\pa_\J)
\end{equation}
For simplicity, we also denote by $\ol\pa$ the operator $\ol\pa_\J$ in this section.
\bgn{definition}\label{a family Fut}
For $u\in C^\infty(M, \C)$, we define $\X^{0,1}_u$ to be
 $\J_\psi\ol\pa u\in \ol\L_\J$.
 Let $\X_u$ denotes a real element 
 $\X_u^{0,1}+\ol{\X_u^{0,1}}\in T_M\oplus T^*_M.$
 Then the real $\X_u=v+\t\in T_M\oplus T^*_M$ gives  a family $\{F^u_t\}$ of $\wtil\Diff(M)$ 
 which is given by 
 $F_t^u:= e^{t d\t }f_{v,t},$ where $f_{v, t}$ is a family of diffeomorphisms generated by the vector field $v$ and $d\t$ is a $d$-exact $2$-form.
 \end{definition}
Let $\{\J_t^u\}$ be deformations of generalized complex structures which are given by
$\J_0^u=\J$ and $\J_t^u :=(F_t^u)_{\#}\J,$
where $(F_t^u)_{\#}\J$ denotes the action of $F^u_t\in \wtil\Diff(M)$ on $\J.$
Then infinitesimal deformations of $\J_t^u$ is given by 
$\e_u:=\ol\pa\J_\psi\ol\pa u\in \w^2\ol\L_\J.$

From the moment map formula, we already know the formula of derivation of generalized scalar curvature under deformations of $\J$ preserving $\J_\psi.$
 For $w\in C^\infty(M, \R)$, we have 
$$
\frac{d}{dt}\int_M S(\J_t^u)\, w\, \vol_M |_{t=0}=\Ome_\B(\dot\J_{\e_u}, \dot\J_{\e_w}),
$$
where $\Ome_\B$ is the K\"ahler form on $\B_{\J_\psi}(M)$ as before and $\dot\J_{\e_u}$ and $ \dot\J_{\e_w}$
are infinitesimal deformations given by $\e_u$ and $\e_w$, respectively.
Then $\Ome_\B(\dot\J_{\e_u}, \dot\J_{\e_w})$ is given by the Imaginary part of the integration of $h(\e_u, \e_w)$ over $M$,
$$
\Ome_\B(\dot\J_{\e_u}, \dot\J_{\e_w})=\,\Im \int_M h(\e_u, \e_w)\vol_M=
\frac1{2\sqrt{-1}}\{\int_M h(\e_u, \e_w)-h(\ol\e_u, \ol \e_w)\}\vol_M,
$$
and $h(\,,\,)$ denotes the Hermitian metric on $\w^2\ol\L_\J$ which is given by 
$$
h(\e_1, \e_2):={4}\,\tr\(\ad_{\e_1},\,\,\ol\ad_{\e_2}\),
$$
where $\ad_{\e_1}:=[\e_1, \,]\in \Hom(\L_\J, \ol\L_\J)$ and $\ol\ad_{\e_2}:=[\ol \e_2, \,\,]\in\Hom(\ol\L_\J, \L_\J)$
(see also Section 7 in \cite{Goto_2016}).
Since 
$h(\e_u, \e_w)={4}\,h(\ol\pa\J_\psi\ol\pa u, \,\, \ol\pa\J_\psi\ol\pa w),$
$h(\ol\e_u,\ol\e_w)={4}\,h(\pa\J_\psi\pa \ol u, \,\,\pa\J_\psi\pa\ol w)$, 
applying (\ref{L=(olpaJpsiolpa}) and (\ref{ol L=(paJpsi}), we obtain
\bgn{align}\label{(-i)intM Luw-ol L}
\frac{d}{dt}\int_M S(\J_t^u)\, w\, \vol_M |_{t=0}=\frac{2}{\sqrt{-1}}\,\{\int_M (Lu)w\, \vol_M -
\int_M (\ol L \ol u)\ol w\, \vol_M\}
\end{align}
A complex function $u$ is written as $u=u_{\Re}+\sqrt{-1}u_{\ss\Im}.$
If $u$ is a real function, i.e., $u=u_{\Re}$,
then we have $\X_u:=\J_\psi\ol\pa u+\J_\psi \pa u=\J_\psi du.$
Thus $\X_u$ is a 
Hamiltonian element of $u$ with respect to $\J_\psi.$ 
We denote by  $\{F^{u_{\Re}}_t\}$
the corresponding family of $\wtil\Diff(M)$ which gives $\J_t^u:=(F_t^{u_{\Re}})_\#\J$. 
Since $(F_t^{u_{\Re}})$ also preserves the volume form, we obtain
\bgn{align}
S(\J_t^u)w\, \vol_M=&(F_t^{u_{\Re}})_\# (S(\J_0)\,\vol_M)w\\
=& (F_{t}^{u_{\Re}})^*(S(\J_0)\,\vol_M)w.
\end{align} 
(Note that $F_\#=F^{-1}_*\oplus F^*$ acts on a differential form $\a$ by $F_\# \a=F^*\a$.)
Since $F_{t}^{u_{\Re}}$ is a Hamiltonian element of $u_{\Re}$, we have
\bgn{align}
\frac{d}{dt}\int_M S(\J_t^u)\, w\, \vol_M \Big|_{t=0}=+\int_M \{ u_{\Re}, \,\,S(J_0)\}_{\J_\psi}w\, \vol_M,
\end{align}
where $\{\,,\,\}_{\J_\psi}$ denotes the Lie derivative
$\L_{\X_u}S(\J_0)$
which is a generalization of
the Poisson bracket.
Thus if both $u$ and $w$ are real functions, i.e.,
$u=u_{\Re}, w=w_{\Re}$, then from (\ref{(-i)intM Luw-ol L}) we have 
\bgn{equation}\label{intMuRe S(J0)poi wRe}
\int_M \{ u_{\Re}, \,\,S(\J_0)\}_{\J_\psi}w_{\Re}\, \vol_M=
\frac2{\sqrt{-1}}\int_M (Lu_{\Re}-\ol L u_{\Re})w_{\Re}\, \vol_M 
\end{equation}
Then we obtain
\bgn{proposition}\label{-iLolL=isJ}
$$\{u, S(\J)\}_{\J_\psi}=\frac{2}{\sqrt{-1}}(L-\ol L)u $$ for every complex function $u$
\end{proposition}
\bgn{proof}
From (\ref{intMuRe S(J0)poi wRe}),
the formula holds for a real function $u$.
Then it follows that the formula holds for every complex function since both sides are $\C$-linear with respect to $u$.
\end{proof}
Next we shall show the derivation formula of the generalized scalar curvature in the case of a pure imaginary function $u=\sqrt{-1}u_{\Im}$.
Then 
$\J_\psi\ol\pa u=\sqrt{-1}\J_\psi\ol\pa u_{\Im}$ gives a family $\{F_t^{u}\}$ of  $\wtil{\Diff}(M)$ 
which yields deformations $\J_t^{u}:=(F_t^{u})_\#\J$ as in Definition \ref{a family Fut}.
Then we have 
\bgn{proposition}\label{Jtu|t=0=L+olL u}
$$\frac{d}{dt}S(\J_t^{u})|_{t=0}=2(L+\ol L)u_{\Im}$$
\end{proposition}
\bgn{proof}
Since $\ol\pa\J_\psi\ol\pa u=2\sqrt{-1}\ol\pa_-\ol\pa_+ u\in \ol\L_\J^+\w\ol\L_\J^-,$ we see that 
$\ol\pa\J_\psi\ol\pa u$ is an infinitesimal tangent of generalized complex structures at $\J$ preserving 
$\J_\psi.$
Thus we can apply the formula of derivation of the moment map.
Since the derivation of the generalized scalar curvature is given by the Moment map formula as before, we have 
$$
\frac{d}{dt}\int_M S(\J_t^{u})\, w\, \vol_M |_{t=0}=\Ome_\B(\dot\J_{\e_u}, \dot\J_{\e_w})=\,\Im \int_M h(\e_u, \e_w)\vol_M,
$$
for a real function $w\in C^\infty(M, \R).$
Since $u=\sqrt{-1}u_{\Im}$, we have 
\bgn{align}h(\e_u, \e_w)=&4h(\sqrt{-1}\ol\pa\J_\psi\ol\pa u_{\Im}, \,\, \ol\pa\J_\psi\ol\pa w),\\
h(\ol\e_u,\ol\e_w)=&4h(-\sqrt{-1}\pa\J_\psi\pa  u_{\Im}, \,\,\pa\J_\psi\pa w).
\end{align}
Thus from (\ref{(-i)intM Luw-ol L}), we also have
\bgn{align}
\frac{d}{dt}\int_M S(\J_t^u)\, w\, \vol_M \Big|_{t=0}=2\int_M (Lu_{\Im})w\, \vol_M +
2\int_M (\ol L  u_{\Im}) w\, \vol_M
\end{align}
Thus we obtain the result.
\end{proof}
Then we obtain
\bgn{theorem}\label{the Lie algebra of the reduced automorphisms }
If the scalar curvature $S(\J, \J_\psi)$ is a constant, then the Lie algebra of the reduced automorphisms 
$\frak g_{\red}$ is reductive.
\end{theorem}
\bgn{proof}
Since $S(\J)$ is a constant, we have $\{u,  S(J)\}_{\J_\psi}=0.$
Then from Proposition \ref{-iLolL=isJ}, we have 
$L=\ol L.$ 
Thus if a complex function $u$ satisfies $Lu=0,$ then $\ol L u=0.$
Then the real part  $u_{\Re}$  of $u$ also satisfies $Lu_{\Re}=0.$
Hence 
if a complex function $u$ satisfies $\ol\pa \J_\psi \ol\pa u =0, $
then both $u_{\Re}$ and $u_{\Im}$ also satisfy 
\bgn{equation}
\qquad \qquad \ol\pa\J_\psi \ol\pa u_{\Re}=0, \quad \ol\pa\J_\psi\ol\pa u_{\Im}=0
\end{equation}
Thus the condition (\ref{uRe uIm}) holds.
Then it follows from Proposition \ref{complexification}
that $\frak g_{\red}$ is reductive.
\end{proof}
\bgn{theorem}\label{main theorem}
Let $(M, \J)$ be a $2n$ dimensional compact generalized complex manifold. 
We assume that $H^{1}(M,\C)=0.$
If $M$ admits a generalized K\"ahler structure $(\J, \J_\psi)$ of symplectic type with constant 
scalar curvature, then 
the Lie algebra $\frak g_\J$ is a reductive Lie algebra.
\end{theorem}
\bgn{proof}
From Corollary \ref{cor hodd=0},
we see $\frak g_\J =\frak g_{\red}.$
Then the result follows from Theorem \ref{the Lie algebra of the reduced automorphisms }
\end{proof}
\bgn{theorem}\label{main theorem Poisson deformations}
Let $(M, I, \ome)$ be a compact K\"ahler manifold with a holomorphic Poisson structure $\b\neq 0.$
We assume $H^{1}(M, \C)=0.$
We denote by $(M, \J_{\b t}, \J_{\psi_t})$  a generalized K\"ahler manifold given by Poisson deformations.
Then if the scalar curvature $S(\J_{\b t}, \J_{\psi_t})$ is a constant, the Lie algebra of the automorphisms $\frak g_{\J_{\b t}}$ is a reductive Lie algebra. 
\end{theorem}
\bgn{proof}
It follows from Theorem \ref{the Lie algebra of the reduced automorphisms 
} that $\frak g_{\red}$ is reductive. 
Then from Proposition \ref{reduced auto coincides with}, we have 
$\frak g_{\red}=\frak g_{\J_{\b t}}.$
Thus the result follows.
\end{proof}

\section{The Lie algebra $\frak g_\b$ of automorphisms of $(\C P^2, \J_\b)$}\label{The Lie alg CP2}
Let $X=(M, J, \ome)$ be the complex projective surface $\C P^2$ and $\b$ a Poisson structure on $X$,
where $M$ denotes the underlying differentiable manifold and $J$ is a complex structure and $\ome$ is a K\"ahler structure on $M.$
Since $\b$ is a holomorphic section of $K_X^{-1}$, 
$\b$ is given by a homogeneous polynomial $f(z_0, z_1, z_2)$ of degree $3$ (a cubic curve).
It is known that cubic curves are classified
into nine cases as shown in the following figures:

\clearpage
%%%%%%%%%%%%%%%%%%%%%%%
\bgn{center}
{\tiny
\tikzset{every picture/.style={line width=0.1pt}} %set default line width to 0.75pt        

\begin{tikzpicture}[x=0.6pt,y=0.6pt,yscale=-0.6,xscale=0.6]
%[x=0.75pt,y=0.75pt,yscale=-1,xscale=1]
%uncomment if require: \path (0,848); %set diagram left start at 0, and has height of 848

%Curve Lines [id:da7510928091548245] 
\draw    (505.95,63) .. controls (508.89,95) and (542.84,127) .. (581.86,119) ;
%Shape: Circle [id:dp7296799727360935] 
\draw   (52.06,403.5) .. controls (52.14,355.18) and (91.37,316) .. (139.7,316) .. controls (188.02,316) and (227.13,355.18) .. (227.06,403.5) .. controls (226.98,451.82) and (187.75,491) .. (139.42,491) .. controls (91.1,491) and (51.99,451.82) .. (52.06,403.5) -- cycle ;
%Shape: Circle [id:dp024791981184383083] 
\draw   (262.06,402.5) .. controls (262.14,354.18) and (301.37,315) .. (349.7,315) .. controls (398.02,315) and (437.13,354.18) .. (437.05,402.5) .. controls (436.98,450.82) and (397.74,490) .. (349.42,490) .. controls (301.1,490) and (261.99,450.82) .. (262.06,402.5) -- cycle ;
%Shape: Circle [id:dp6547581608096203] 
\draw   (445.06,126.5) .. controls (445.14,78.18) and (484.37,39) .. (532.69,39) .. controls (581.02,39) and (620.13,78.18) .. (620.05,126.5) .. controls (619.98,174.82) and (580.74,214) .. (532.42,214) .. controls (484.1,214) and (444.98,174.82) .. (445.06,126.5) -- cycle ;
%Curve Lines [id:da3386684098800756] 
\draw    (298.54,128) .. controls (299.58,101) and (325.72,17) .. (375.48,166) ;
%Curve Lines [id:da003604874804210012] 
\draw    (298.54,128) .. controls (299.48,165) and (331.48,169) .. (376.65,63) ;
%Curve Lines [id:da1595220823813026] 
\draw    (509.77,175) .. controls (537.83,132) and (533.84,127) .. (581.86,119) ;
%Shape: Ellipse [id:dp3748153381225241] 
\draw   (67,123.46) .. controls (67,74.83) and (106.18,35.4) .. (154.5,35.4) .. controls (202.82,35.4) and (242,74.83) .. (242,123.46) .. controls (242,172.1) and (202.82,211.52) .. (154.5,211.52) .. controls (106.18,211.52) and (67,172.1) .. (67,123.46) -- cycle ;
%Shape: Ellipse [id:dp8023559394046058] 
\draw   (256,123.47) .. controls (256,74.83) and (295.18,35.41) .. (343.5,35.41) .. controls (391.82,35.41) and (431,74.83) .. (431,123.47) .. controls (431,172.1) and (391.82,211.53) .. (343.5,211.53) .. controls (295.18,211.53) and (256,172.1) .. (256,123.47) -- cycle ;
%Shape: Ellipse [id:dp21490828422959762] 
\draw   (469,406.48) .. controls (469,357.84) and (508.18,318.41) .. (556.5,318.41) .. controls (604.82,318.41) and (644,357.84) .. (644,406.48) .. controls (644,455.11) and (604.82,494.54) .. (556.5,494.54) .. controls (508.18,494.54) and (469,455.11) .. (469,406.48) -- cycle ;
%Curve Lines [id:da7948762689968605] 
\draw    (98,119.93) .. controls (158.8,64.58) and (149.8,182.33) .. (198,119.93) ;
%Shape: Circle [id:dp8751788903922799] 
\draw   (50.06,667.5) .. controls (50.14,619.18) and (89.37,580) .. (137.7,580) .. controls (186.02,580) and (225.13,619.18) .. (225.06,667.5) .. controls (224.98,715.82) and (185.75,755) .. (137.42,755) .. controls (89.1,755) and (49.99,715.82) .. (50.06,667.5) -- cycle ;
%Shape: Circle [id:dp48923050255461664] 
\draw   (260.06,666.5) .. controls (260.14,618.18) and (299.37,579) .. (347.7,579) .. controls (396.02,579) and (435.13,618.18) .. (435.05,666.5) .. controls (434.98,714.82) and (395.74,754) .. (347.42,754) .. controls (299.1,754) and (259.99,714.82) .. (260.06,666.5) -- cycle ;
%Shape: Ellipse [id:dp44689155445937945] 
\draw   (467,670.48) .. controls (467,621.84) and (506.18,582.41) .. (554.5,582.41) .. controls (602.82,582.41) and (642,621.84) .. (642,670.48) .. controls (642,719.11) and (602.82,758.54) .. (554.5,758.54) .. controls (506.18,758.54) and (467,719.11) .. (467,670.48) -- cycle ;
%Shape: Ellipse [id:dp15455345250424601] 
\draw   (104.56,403.5) .. controls (104.56,388.86) and (120.23,377) .. (139.56,377) .. controls (158.89,377) and (174.56,388.86) .. (174.56,403.5) .. controls (174.56,418.14) and (158.89,430) .. (139.56,430) .. controls (120.23,430) and (104.56,418.14) .. (104.56,403.5) -- cycle ;
%Straight Lines [id:da6480909600528805] 
\draw    (89.56,353.5) -- (189.56,453.5) ;
%Shape: Ellipse [id:dp16183392190678936] 
\draw   (314.56,402.5) .. controls (314.56,391.45) and (330.23,382.5) .. (349.56,382.5) .. controls (368.89,382.5) and (384.56,391.45) .. (384.56,402.5) .. controls (384.56,413.55) and (368.89,422.5) .. (349.56,422.5) .. controls (330.23,422.5) and (314.56,413.55) .. (314.56,402.5) -- cycle ;
%Straight Lines [id:da939034334053925] 
\draw    (280,372.6) -- (380,472.6) ;
%Straight Lines [id:da8021902951311335] 
\draw    (523,343.6) -- (623,443.6) ;
%Straight Lines [id:da04979720982457447] 
\draw    (505.4,423.8) -- (625.4,405.8) ;
%Straight Lines [id:da5928323376598936] 
\draw    (555.4,336.8) -- (534.4,472.8) ;
%Straight Lines [id:da8866397393465051] 
\draw    (89.56,626.5) -- (189.56,726.5) ;
%Straight Lines [id:da09338936490963445] 
\draw    (73.8,666.4) -- (197.8,677.4) ;
%Straight Lines [id:da9622429217687143] 
\draw    (85.8,715.4) -- (178,630.2) ;
%Straight Lines [id:da6601303821432841] 
\draw    (292.76,667.88) -- (417.76,666.88) ;
%Straight Lines [id:da45589536483094706] 
\draw    (287.76,671.88) -- (419.76,671.88) ;
%Straight Lines [id:da03358603236529323] 
\draw    (297.56,616.5) -- (397.56,716.5) ;
%Straight Lines [id:da2207398314576795] 
\draw    (488.38,658.78) -- (624.38,658.78) ;
%Straight Lines [id:da8117725580453774] 
\draw    (488.62,670.38) -- (620.38,670.58) ;
%Straight Lines [id:da2725423069472146] 
\draw    (492,664.68) -- (620.76,664.88) ;

% Text Node
\draw (71.91,496.9) node [anchor=north west][inner sep=0.75pt]   [align=left] { \ \ \  Figure 4.\\ conic and line \\in general position\\};
% Text Node
\draw (473.9,227.9) node [anchor=north west][inner sep=0.75pt]   [align=left] { \ \ \ Figure 3.\\ \ \ \  cuspidal curve\\};
%(473.9,227.9)
% Text Node
\draw (289.91,500.9) node [anchor=north west][inner sep=0.75pt]   [align=left] { \ \ \  Figure 5.\\ \ \  conic and line \\  in special position\\};
% Text Node
\draw (76,226.66) node [anchor=north west][inner sep=0.75pt]   [align=left] { \ \ \ Figure 1. \\\ \ \ smooth  \\elliptic curve\\};
% Text Node
\draw (488,505.7) node [anchor=north west][inner sep=0.75pt]   [align=left] { \ \ \  Figure 6.\\ \ \  three lines \\ \ in general position\\};
% Text Node
\draw (279,231.66) node [anchor=north west][inner sep=0.75pt]   [align=left] { \ \ \  Figure 2.\\ \ \   nodal curve\\};
% Text Node
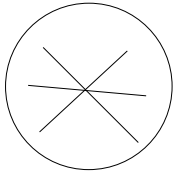
\draw (62.91,766.9) node [anchor=north west][inner sep=0.75pt]   [align=left] { \ \ \ Figure 7.\\ three lines intersect\\ \ \ \ at one point \\};
% Text Node
\draw (259.91,769.9) node [anchor=north west][inner sep=0.75pt]   [align=left] { \ \ \ \ \ Figure 8.\\ \ \ \ the double line \\ \ \ \ and a line\\};
% Text Node
\draw (486,769.7) node [anchor=north west][inner sep=0.75pt]   [align=left] { \ \ \  Figure 9.\\ \  the triple line\\};

\end{tikzpicture}

}\end{center}
From the deformations-stability theorem, there exists a family of generalized K\"ahler structures $(\J_{\b t}, \J_{\psi_t}).$
Since $H^1(X,\C)=0,$ it follows from Proposition \ref{reduced auto coincides with}  that 
the Lie algebra of automorphisms $\frak g_{\J_{\b t}}$ 
coincides with the Lie algebra of reduced automorphisms 
$\frak g_{\red}$.
Since $\frak g_{\J_{\b t}}$ is the same for  $t\neq 0$, for simplicity, we denote by $\frak g_\b$ the Lie algebra $\frak g_{\J_{\b t}}$ for $t\neq 0$.
From Theorem \ref{main theorem}, if the scalar curvature $S(\J_{\b t}, \J_{\psi_t})$ is constant, the Lie algebra $\frak g_{\b}$ is reductive. 
Thus we have an obstruction to the existence of constant scalar curvature on a generalized complex manifold $(M, \J_{\b t})$ for $t\neq 0$.
From Proposition \ref {Poisson vector fields}, 
we already know that  the Lie algebra ${\frak g}_{{\b }}$ is given by the Lie algebra of holomorphic vector fields 
preserving the holomorphic Poisson structure $\b$, i. e., ${\frak g}_{{\b }}=\{ V\in H^0(X, T^{1,0}_J)\, |\, \L_V\b=0\, \}.$
The Lie algebra {\it sl}$(3, \C)$ acts on $\C P^2$ linearly. 
Let  $f:=f(z_0, z_1, z_2)$ by a homogeneous polynomial of degree $3$ given by $\b$.
Then $\frak g_\b$ is given by 
$$
\frak g_\b =\{\, a \in \text{\it sl}(3, \C)\, |\, 
a^*f=0\, \},
$$
where $a^*f$ denotes the action of the Lie algebra sl$(3, \C)$ on the space of homogeneous polynomials of degree $3$.

Then we have the following explicit calculations: \\ \\
{\bf  Figure 1}:
If an anticanonical divisor is a smooth elliptic curve $C$, then we see $\frak g_\b=0.$ \\
{\bf Figure 2}:
In the case of a nodal curve, we also have $\frak g_\b=0$. 
\\
{\bf Figure 3,4,5,6}\label{three lines}:
If an anticanonical divisor is a curve given in Figure 3, 4, 5, 6, then it follows that $\frak g_\b$
is abelian. Thus $\frak g_\b$ is reductive.
\\
{\bf Figure 7}\label{three lines at one point}:
In the cases where  an anticanonical divisor $f$ is given by three lines intersecting at one point (Figure 7), 
$f$ can be taken a following form $z_0 z_1(z_0+z_1).$
Then $\frak g_{\b}$ is generated by matrices 
 $$
\bgn{pmatrix}0&0&0\\ 0&0&0\\ *&*&0\end{pmatrix}
$$
Then $\frak g_\b$ has a nonabelian solvable ideal. Thus $\frak g_\b$ is not reductive.
\\
{\bf Figure 8}\label{double line and line}:
If the anticanonical divisor $f$ is $z_0^2z_1$ which consists of a double line $z_0^2$ and another line $z_1$ (Figure 8). 
Then we see $\frak g_\b$ is generated  by 
$$
\bgn{pmatrix}-2t&0&0\\ 0&t&0\\ *&*&t\end{pmatrix} , \qquad \text{\rm  } (t\in \C).
$$
Thus $\frak g_\b$ is not reductive also.
\\{\bf Figure 9}\label{triple line}:
If an anticanonical divisor $f$ is a triple line $z_0^3$ (figure. 9),
then $\frak g_\b$ is generated by the following elements
$$
\bgn{pmatrix}0&0&0\\ *&t&*\\ *&*&-t\end{pmatrix}, \qquad (t\in \C).
$$
Then $\frak g_\b$ has a nonabelian solvable ideal which is generated by 
$$
\bgn{pmatrix}0&0&0\\ *&0&0\\ *&0&0\end{pmatrix}
$$
Thus the radical of $\frak g_\b$ is not abelian. 
Thus $\frak g_\b$ is not reductive.
\\
Thus we see that $\frak g_\b$ is reductive for the cases as in Figures\,1, 2, 3, 4, 5, 6.
However $\frak g_\b$ is not reductive for the cases as in figures\,7 ,8, 9. 
Then we have
\bgn{proposition}
Let $(M, \J_\b)$ be a generalized complex manifold which $\b$ is given by  three cases as in Figures\,7,8,9.
Then $(M,\J_\b)$ does not admits a generalized K\"ahler structure with constant scalar curvature. 
\end{proposition}
\bgn{proof}
Since $\frak g_\b$ is not reductive in these three cases, 
the result follows Theorem \ref{main theorem}.
\end{proof}
We have previously shown that the existence of generalized K\"ahler structures with constant scalar curvature in the cases of Figure 6.
\bgn{proposition}
If $\b$ is given by three lines in general position, then there exists generalized K\"ahler structures with constant scalar curvature 
\end{proposition}
\bgn{proof}
In the cases when the Poisson tensor is given by an action of $2$-dimensional torus on $\C P^2$, we can apply the result in \cite{Goto_2016} (see Proposition 12.3).
\end{proof}
\section{Deformations of generalized K\"ahler manifolds with constant scalar curvature}
\label{DeformcscGK}
Let $(\J, \J_\psi)$ be a generalized K\"ahler manifold of symplectic type on a compact manifold $M$ and $S(\J, \J_\psi)$ the scalar curvature of $(\J, \J_\psi).$
For simplicity,  $(\J, \J_\psi)$ is denoted by $(\J, \psi).$
We assume that the scalar curvature $S(\J, \psi)$ is a constant $\h S$, which is topologically given by the $1$-st Chern class of 
the canonical line bundle $K_\J$ together with the class $[\psi].$
We will consider a natural deformation problem for generalized K\"ahler structures of symplectic type with constant scalar curvature under fixing the cohomology class $[\psi].$
If we fix the class $[\psi]$ and $\J$ and deform $\psi$ such that $(\J, {\psi_s})$ are generalized K\"ahler structures parametrized by $s\in (-\e, \e)$, then it follows from the $\pa\ol\pa$-lemma for generalized K\"ahler manifolds that the derivative $\dot{\psi}:=\frac{d}{dt}\psi_s|_{s=0}$ is given by 
\bgn{equation}\label{dotpsi:=dolpaucdot}
\dot\psi:=d(\ol\pa_\J u\cdot\psi) \in U^{0,-n}\oplus U^{0, -n+2},
\end{equation}
for a complex function $u.$
However, we only consider deformations $\psi_s$ given by a real function $u$ to apply the implicit function theorem later.  
If we have deformations $\{\J_t\}$ of generalized complex structures with $\J_0=\J$ and $t$ denotes a parameter of deformations satisfying
$|t|<\e$.
Then deformation-stability theorem provides  deformations of generalized K\"ahler structures $(\J_t, \J_{\psi_t})$
such that $\J_{\psi_0}=\J_\psi.$ 
Further it turns out that deformation-stability theorem yields 2-parameter deformations 
$(\J_t, {\psi_{t, u}})$ of generalized K\"ahler manifolds of symplectic type which are smoothly parametrized by $t$ and a real function $u$. 
We denote by $S(\J_t, \psi_{t, u})$ the scalar curvature of $(\J_t, \psi_{t, u}).$
We normalize a function $u$ such that $\int_M u \vol_M=0.$
Let $L_k^2(M)$ be the Sobolev space of real functions on $M$ whose 
first $k$ derivatives are square integrable. The Sobolev embedding theorem states that 
$L_k^2(M)\subset C^l(M)$ if $k>n+l,$ where $2n=\dim_\R M$ and $C^l(M)$ denotes 
the space of continuous functions whose derivative of order at most $l$ are also continuous.
Note that $L_k^2(M)$ is a Banach algebra if $k>n.$
We denote by $L^2_k(M)/\R$ the space of normalized functions of $L^2_k(M)$.
We shall define the scalar curvature operator ${\mathcal S}$ as a non-linear differential operator.
First assume that there exist deformations $\{\J_t\}$ of generalized complex structures parametrized by $t$ in a neighborhood of the origin of $\C.$
We need to take $(t, u)$
in a small open set  $D_\e\times {\mathcal U}$  of the origin of $\C\times L^2_{k+4}(M)/\R$ 
such that $(\J_t, \psi_{t, u})$ is a generalized K\"ahler structure, where $D_\e=\{t\in \C\, |\, -\e<|t|<\e\,\}.$
Then we define the map 
$$
{\mathcal S} : D_\e\times  \mathcal U \to L^2_k(M)/\R
$$
which is given by $\mathcal S(\J_t, \psi_{t, u}):= S(\J_t, \psi_{t,u})-\h S.$
Then we have 
\bgn{theorem}\label{implicit function theorem}
For $k>n,$ the map ${\mathcal S}$ is well-defined and the derivative of $\mathcal S$ at the origin
along the direction of the function space $\mathcal U$
is given by 
\bgn{align}
d\mathcal S_{(0,0)}: L^2_{k+4}(M)/\R&\to L^2_k(M)/\R\\
u \qquad&\mapsto 2L u,
\end{align}
where $L$ is the fourth order differential operator 
$L=(\ol\pa_\J\J_\psi\ol\pa_\J)^*\circ(\ol\pa_\J\J_\psi\ol\pa_\J)$
as in (\ref{L=(olpaJpsiolpa}).
\end{theorem}
\bgn{proof} Since $(\ol\pa_\J u\cdot\psi)=(\ol\pa_+ u) \cdot\psi=-\sqrt{-1}(\J_\psi\ol\pa_\J u)\cdot\psi$
, from (\ref{dotpsi:=dolpaucdot}), we have 
$$\dot{\psi}=d(\ol\pa_\J u\cdot\psi)=-\sqrt{-1}d(\J_\psi\ol\pa_J u)\cdot\psi.$$
We denote by $v$ the pure imaginary function $-\frac{\sqrt{-1}}2u.$
Let $F_s^{v}$ be a family of $\wtil\Diff(M)$ for the pure imaginary function $v$ as in
Definition \ref{a family Fut}, where $s$ denotes a parameter. 
Define $\psi_s$ by $(F^v_s)^*\psi.$ 
Since $\frac{d}{ds}(F^v_s)^*(\psi)|_{s=0}$ is given by the Lie derivative $\{d, \,\, (\J_\psi \ol\pa_\J v+\J_\psi\pa_\J \ol v)\}\psi$,
then it follows that 
\bgn{align}
\frac{d}{ds}(F^v_s)^*(\psi)|_{s=0}
=&d(\J_\psi\ol\pa_\J v-\J_\psi \pa_\J v)\cdot\psi=-\sqrt{-1}d (\J_\psi \ol\pa_\J u)\cdot\psi.
\end{align}
Thus we have $$\dot\psi :=\frac{d}{ds}(F^v_s)^*(\psi)|_{s=0}$$ and we have
$\frac{d}{ds}S(\J_t, \psi_{t, us})|_{s=0}=\frac{d}{ds}S(\J_t, (F^v_s)^*\psi)|_{s=0}.$
Since the scalar curvature is equivalent under the action of $\wtil\Diff(M),$
we have 
\bgn{align}
\mathcal S\((F^v_s)_{\#}\circ(F^{-v}_s)_{\#}\J_t, \,\,(F^v_s)^*\psi\)
=&(F_s^u)^*\mathcal S ((F^{-v}_s)_{\#}\J_t, \,\,\psi)
\end{align}
From our assumption $\mathcal S(\J, \psi)=0$ , then we have 
\bgn{align}
\frac{d}{ds} \mathcal S(\J_t, \psi_{t, us})|_{s, t=0}=&\frac{d}{ds}
\mathcal S ((F^{-v}_s)_{\#}\J_t, \,\,\psi)|_{s, t=0}
\end{align}
Since $S(\J, \psi)$ is a constant,
Proposition \ref{-iLolL=isJ} shows $L=\ol L.$
Applying Proposition \ref{Jtu|t=0=L+olL u} to a pure imaginary function $-v=\frac{\sqrt{-1}}2u$, we obtain 
$$\frac{d}{ds}\mathcal S ((F^{-v}_s)_{\#}\J_t, \,\,\psi)|_{s=0}=(L+\ol L)u =2L u$$

Thus the differential of $\mathcal S$ at $(0,0)$ is given by  $d\mathcal S_{0,0}(u)=2Lu.$

\end{proof}
\bgn{theorem}\label{existenceresult}
Let $(\J, \psi)$ be a generalized K\"ahler structure of symplectic type on a compact manifold $M$
with constant scalar curvature $S(\J, \psi).$
We assume that the Lie algebra of the reduced automorphisms $\frak g_{\red}$ is trivial.
Then for deformations generalized complex structures $\{\J_t\}$, ($-\e<|t|<\e),$
there exist deformations of generalized K\"ahler structures $(\J_t, \psi_{t, u_t})$ with constant scalar curvature for sufficiently small $t.$
\end{theorem}
\bgn{proof} From Theorem \ref{implicit function theorem}, $\ker d\mathcal S_{(0,0)}$ is 
given by $\ker L.$
Since $\frak g_{\red}=0,$ it follows $\ker L=0.$
Since $L$ is a self-dual operator, it follows that $d\mathcal S_{(0,0)}$ is surjective and has a right inverse operator. Thus applying the implicit function theorem of Hilbert spaces, we obtain deformations of 
generalized K\"ahler structures $\{(\J_t, \psi_{u_t})\}$ such that $\mathcal S(\J_t, \psi_{u_t})=0.$
Hence we obtain the result. 
\end{proof}
\bgn{example}
Let $(M, J, \ome)$ be a compact K\"ahler manifold with constant scalar curvature. 
We assume that the reduced automorphisms of $(M,J)$ is trivial. 
If there exists a nonzero holomorphic Poisson structure $\b$ on $(M,J),$
then there exist deformations of generalized K\"ahler structures $(\J_{\b t}, \psi_t)$
with constant scalar curvature. In particular, del Pezzo surfaces with trivial automorphisms admit 
generalized K\"ahler structures with constant scalar curvature.
\end{example}

\section{Generalized extremal K\"ahler manifolds and Hessian formula}
\label{GEXKandHessian}
Our moment map framework naturally leads us to 
define a generalized extremal K\"ahler manifold.
We denote by $\J_\psi$ a generalized complex structure given by $\psi$ which is fixed in this section.
\bgn{definition}
Let $(M, \J, \J_\psi)$ be a generalized K\"ahler manifold and 
$S(\J)$ the generalized scalar curvature of $(M, \J, \J_\psi).$
If $\ol\pa\J_\psi \ol\pa S(\J)=0,$ then $(\J, \J_\psi)$ is 
a {\it generalized extremal K\"ahler structure}
and $(M, \J, \J_\psi)$ is called {\it a generalized extremal K\"ahler manifold.}
\end{definition}
\bgn{proposition}
 Let $\E$ be a functional on the space $$\{\J \, |\, (\J, \J_\psi) \,\,\text{\rm  is a generalized K\"ahler structure}\}$$ which is defined by 
 $$
 \E(\J)=\frac12\int_M S(\J)^2 \vol_M$$
 Then a critical point of the functional $\E$ is attained by a generalized extremal K\"ahler structure. 
\end{proposition}
\bgn{proof}
Let $\e\in \w^2\ol\L_\J$ is an infinitesimal deformation which is given by a one-parameter deformation
$\{\J_t\}.$
Then
the differential of $\E$ is given by 
$$
\frac{d}{dt}\E(\J_t)\Big|_{t=0}=\int_M \dot S(\J)S(\J)\vol_M
$$
Since $S(J)$ gives the moment map, we have 
$$
\frac{d}{dt}\E(\J_t)\Big|_{t=0}=\Ome_\B(\dot\J_{\e_{S(\J)}}, \dot\J_\e),
$$
where $\dot\J_{\e_{S(\J)}}$ denotes the infinitesimal deformation corresponding to 
$\ol\pa\J_\psi\ol\pa S(J)$. Since $\Ome_\B$ is nondegenerate, $\J$ is a critical point of $\E$ if and only if 
$\E_{S(\J)}=\ol\pa\J_\psi\ol\pa S(\J)=0.$
Thus the result follows.
\end{proof}

In order to calculate the Hessian of the functional $\E$, we need several lemmas. A diffeomorphism 
$f\in \Diff(M)$ gives rise to a bundle map $f_\#: \TT\to \TT$ which is defined by
$f_\#(v, \t) =(f_*^{-1}(v), f^*(\t)),$
where $v\in T_M$ and $\t\in T^*_M$. 
A $2$-form $b$ acts on $\TT$ by $\Ad_{e^b}.$
Then 
$F=e^{b}f\in \wtil\Diff(M)$ also gives a bundle map which is the composition
$F_\#(v, \t)=\Ad_{e^b}\circ f_\#(v, \t)=\Ad_{e^b}(f_*^{-1}(v), f^*(\t))=(f_*^{-1}(v)+b(f^{-1}_*(v)), f^*(\t))$.
Thus $F\in \wtil\Diff(M)$ acts on $\J$ by the adjoint $F_\#\circ \J\circ F_\#^{-1}.$
For simplicity we denote by $F_\#\J$ the adjoint $F_\#\circ \J\circ F_\#^{-1}.$
\bgn{lemma}\label{FL(J, Jpsi)}
We denote by $P:=P(\J, \J_\psi)$ the operar $\ol\pa_\J \J_\psi\ol\pa_\J$ and 
$L:=L(\J, \J_\psi)$ the $4$-th order differential operator $(\ol\pa_\J \J_\psi\ol\pa_\J)^*(\ol\pa_\J \J_\psi\ol\pa_\J).$
For $F=e^{d\eta}f\in \wtil\Diff(M),$
we have 
$$
P(F_\#\J, F_\#\J_\psi)=F_\#\circ P(\J, \J_\psi)\circ F_\#^{-1}
$$
Further we have 
$$L(F_\#\J, \,\,F_\#\J_\psi)=F_\#\circ L(\J, \J_\psi)\circ F_\#^{-1}
$$
\end{lemma}
\bgn{proof}
$F\in \wtil\Diff(M)$ acts on a differential form 
$\a$ by $F_\#\a:=e^b\w f^*\a$, where $f^*$ denotes the pull back of $\a$ by $f\in \Diff(M)$ and 
$e^b\w$ is the wedge product of $e^b$.
The we have $d\circ F_\#=F_\#\circ d$.
Since  $F_\#\J=F_\#\circ \J\circ F_\#^{-1}$,  we have 
 \bgn{align*}2\ol\pa_{F_\#(\J)}=d+\sqrt{-1}F_\#(\J) d=&
F_\#\circ d\circ F_\#^{-1}+\sqrt{-1}F_\#\circ \J \circ d\circ F_\#^{-1}\\
=&2F_\#\circ\(\ol\pa_\J\)\circ F_\#^{-1}.
\end{align*}
Since $\F_\#(\J_\psi)=F_\#\circ \J_\psi\circ F_\#^{-1}$, we also have
\bgn{align}
\(\ol\pa_{F_\#\J}\circ F_\#\J_\psi\circ \ol\pa_{F_\#\J} \) =&
\(F_\#\circ \ol\pa_\J\circ \J_\psi\circ \ol\pa_\J\circ F_\#^{-1}\)\\
=&F_\#\(\ol\pa_\J\circ\J_\psi\circ\ol\pa_\J\)F_\#^{-1}
\end{align}
Hence we obtain 
$$
P(F_\#\J, F_\#\J_\psi)=F_\#\circ P(\J, \J_\psi)\circ F_\#^{-1}
$$
We also have 
$(F_\#\circ P\circ F_\#^{-1})^{*_F}=F_\#\circ P^*\circ F_\#^{-1}$.
Then we obtain
\bgn{align}
L({F}_\#(\J), {F}_\#(\J_{\psi}))=&(F_\#\circ P^*\circ {F}_\#^{-1})\circ (F_\# \circ P\circ F_\#^{-1})\\
=&F_\#\circ P^* P\circ F_\#^{-1}\\
=&F_\#\circ L\circ F_\#^{-1}
\end{align}
Hence 
$$
L({F}_\#(\J), {F}_\#(\J_{\psi}))(F_\# S) ={F}_\#( L(\J, \J_\psi)S)
$$
Thus we obtain 
$$L(F_\#\J, \,\,F_\#\J_\psi)=F_\#\circ L(\J, \J_\psi)\circ F_\#^{-1}
$$
\end{proof}
Let $(M, \J, \J_\psi)$ be a generalized extremal K\"ahler manifold as before.
For a real function $u,$ we define a real element $e_u\in \TT$ by 
\bgn{equation}\label{ducdotpsi=}(du)\cdot\psi ={}\sqrt{-1}e_u\cdot\psi.
\end{equation}
Applying $\J$ to the both sides, we see that (\ref{ducdotpsi=}) is equivalent to
 $\J_\psi(du) =-e_u$.

Let $\{ F_s\}$ be the family of the extended diffeomorphisms $\wtil\Diff(M)$ which is generated by $ -\J e_u \in T_M\oplus T_M^*,$ that is, $\frac{d}{ds}F_s|_{t=0}=-\J e_u,$ and $F_0=\id$ and $\e<s<\e$ for a sufficiently small $\e>0.$
Since $\J\J_\psi=\J_\psi\J,$ we have 
$-\J e_u=\J \J_\psi du =\J_\psi \J du =\J_\psi(\sqrt{-1}\,\ol\pa u-\sqrt{-1}\,\pa u).$
Thus $F_s$ coincides with $F_s^{\sqrt{-1}u}$ which are deformations given by a pure imaginary function $\sqrt{-1}u$
as in Definition \ref{a family Fut}.
Since $-\J e_u$ is written as $v+\eta\in T_M\oplus T_M^*$, $F_s$ is given by $F_s=e^{d\eta_s}f_s,$
where $d\eta_s$ is the exact $2$-form and $f_s$ is a diffeomorphism of $M.$
Then $F_s$ acts on $\psi $ by $(F_s)_\#\psi:=e^{d\eta_{s}}\w f_{s}^*\psi.$
We also denote by $\psi_s$ the nondegenerate, pure spinor $e^{d\eta_{s}}\w f_{s}^*\psi.$
For simplicity, we denote by $S$ the scalar curvature $S(\J)$.

%%%%%%%%%%%%%%%%%%%%%%%%%%%%%%%%%%%%%%%%%%%%%%

\bgn{lemma}\label{The projection of LJeupsi}
Let $\{\psi_{-s}\}$ be a family which is given by $(F_{-s})_\#\psi$. 
Then we have 
$$
\frac{d}{ds}\psi_{-s}\big|_{s=0}=\L_{\J( e_u)}\psi=d(\ol\pa_+u-\pa_-u )\psi.
$$
The projection of $\L_{\J( e_u)}\psi$ to the component $U_{\J_\psi}^{-n+2}$
is given by 
$$
\pi_{U^{-n+2}_{\J_\psi}}\L_{\J(e_u)}\psi=-2(\ol\pa_+\pa_-u)\psi.
$$
\end{lemma}
\bgn{proof}
From (\ref{ducdotpsi=}), we have
$$d\J(du)\cdot\psi=+\sqrt{-1}d\J(e_u)\cdot\psi=+\sqrt{-1}\L_{\J(e_u)}\psi.
$$
where $\L_{\J(e_u)}$ denotes the differential operator $d\circ \J(e_u)+\J(e_u)\circ d.$
Then we have
\bgn{align}
\L_{ \J(e_u)}\psi=&-{\sqrt{-1}}d(\J d u)\psi\\
=&-{\sqrt{-1}}d(\sqrt{-1}\,\ol\pa_+ u -\sqrt{-1}\pa_- u)\cdot\psi\\
=&d(\ol\pa_+u-\pa_-u )\psi\\
\end{align}
Then we have 
\bgn{align}
\pi_{U^{-n+2}_{\J_\psi}}\L_{\J(e_u)}\psi=&\pi_{U^{-n+2}_{\J_\psi}}\(d(\ol\pa_+u-\pa_-u )\psi\)\\
=&2(\pa_-\ol\pa_+u)\psi=-2(\ol\pa_+\pa_-u)\psi
\end{align}
\end{proof}
Note that $(\ol\pa_+\pa_-u)\in \ol\L_\J^+\cap \L_\J^-\subset \w^2\ol\L_{\J_\psi}$ which is not a differential operator but a tensor.
Let $\{\J_{\psi_{-s}}\}$ be a family of generalized complex structures which are given by 
$\psi_{-s}.$
We denote by $\dot{\J}_\psi$ the differential $\frac{d}{ds}\J_{\psi_{-s}}\big|_{s=0}$.
\bgn{lemma} 
Let $E_{-s}\in \L_{\psi_{-s}}$ be a smooth family with $E=E_0\in \L_\psi.$
We denote by $\dot{E}$ the differential $\frac{d}{ds}E_{-s}|_{s=0}$ and $\pi_{\L_\psi}$ the projection to 
$\L_\psi.$
Then we have
\bgn{align}\label{dotLpsiE=2}
\dot{\J_\psi }E=&-2\sqrt{-1} \pi_{\ol\L_\psi}\dot E=+4\sqrt{-1}[\ol\pa_+\pa_- u, \,\,E], \qquad \text{\rm for } E\in \L_{\J_\psi},\\
\dot{\J_\psi }\ol E=&+2\sqrt{-1} \pi_{\L_\psi}\ol{\dot E}=-4\sqrt{-1}[\pa_+\ol\pa_- u, \,\,\ol E],
\qquad \text{\rm for }\ol E\in \ol\L_{\J_\psi}.
\end{align}
\end{lemma}
\bgn{proof}
Since $E_{-s}\in \L_{\psi_{-s}}$, we have 
$E_{-s}\cdot \psi_{-s}=0$
Since $E\in\L_\psi$, we also have
$E\cdot\psi=0$. Thus we have
$E\cdot\pi_{U^{-n}_{\J_\psi}}\(\L_{\J (e_u)}\psi\)=0$. Since $\dot{\psi}=\L_{\J (e_u)}\psi$, 
it follows from Lemma \ref{The projection of LJeupsi} that the differential of both sides $E_{-s}\cdot \psi_{-s}=0$ gives 
\bgn{align}
0=&\dot E\cdot\psi +E\cdot\dot \psi=\dot E\cdot\psi +E\cdot\L_{\J (e_u)}\psi\\
=&\dot E\cdot\psi +E\cdot\pi_{U^{-n+2}_{\J_\psi}}\(\L_{\J( e_u)}\psi\)\\
=&\dot E\cdot\psi -2[E, \,\, \ol\pa_+\pa_-u]\cdot\psi\\
=&\dot E\cdot\psi +2[ \ol\pa_+\pa_-u,\,\, E]\cdot\psi,
\end{align}
where $[\,,\,]$ denotes the commutator of the Clifford algebra.
Hence we have
\bgn{equation}\dot E\cdot\psi=-2[ \ol\pa_+\pa_-u,\,\, E]\cdot \psi
\end{equation} 
The $\ol\L_\psi$-component of 
$\dot E$ is denoted by $\pi_{\ol\L_\psi}\dot E$. 
Then we have  
\bgn{equation}\label{piolLpsidotE=-2}\pi_{\ol\L_\psi}\dot E=-2[ \ol\pa_+\pa_-u,\,\, E].\end{equation}
Since $E_{-s}\in \L_{\J_{\psi_{-s}}},$ we have 
$$
\J_{\psi_{-s}}E_{-s}=-\sqrt{-1}E_{-s}.
$$
The differential of the both sides yields 
$$
\dot{\J_\psi}E+\J_\psi\dot E=-\sqrt{-1}\dot E$$
Then from (\ref{piolLpsidotE=-2}), we have

\bgn{equation}\label{dotLpsiE=2}
\dot{\J_\psi }E=-2\sqrt{-1} \pi_{\ol\L_\psi}\dot E=+4\sqrt{-1}[\ol\pa_+\pa_- u, \,\,E]
\end{equation}
Taking the complex conjugate,  we also have 
\bgn{equation}\label{dotJpsiolE=-2}
\dot{\J_\psi }\ol E=2\sqrt{-1} \ol{\dot E}=-4\sqrt{-1}[\pa_+\ol\pa_- u, \,\,\ol E],
\end{equation}
for $\ol E\in \ol\L_\psi$.
\end{proof}
Then applying (\ref{dotLpsiE=2}) and (\ref{dotJpsiolE=-2}) to $\ol\pa_\pm S$ respectively, we have 
\bgn{align}\label{dotJpsiolpaS=}
\dot{\J_\psi}(\ol\pa S)=&\dot{\J_\psi}(\ol\pa_+ S+\ol\pa_- S)\notag\\
=&-4\sqrt{-1}[\pa_+\ol\pa_- u, \,\,\ol\pa_+ S]+4\sqrt{-1}[\ol\pa_+\pa_- u, \,\, \ol\pa_- S]\\
=&+4\sqrt{-1}[\ol\pa_-\pa_+ u, \,\,\ol\pa_+ S]+4\sqrt{-1}[\ol\pa_+\pa_- u, \,\, \ol\pa_- S]\notag
\end{align}
\bgn{lemma}\label{2olpalanpa-uolpa}
If $\ol\pa_+\ol\pa_- S=0$, then we have
\bgn{align}
2\ol\pa_+\lan \pa_- u, \,\, \ol\pa_- S\ran_{\tt}=&[\ol\pa_+\pa_- u, \,\, \ol\pa_- S]\label{20lpalanpa-uolpa-S}\\
2\ol\pa_-\lan \pa_+u, \,\, \ol\pa_+ S\ran_{\tt}=&[\ol\pa_-\pa_+ u, \,\, \ol\pa_+ S]
\end{align}
\end{lemma}
\bgn{proof}
From $
2\lan \pa_- u, \,\, \ol\pa_- S\ran_{\tt}=(\pa_- u)\cdot( \ol\pa_- S)+(\ol\pa_- S)\cdot(\pa_- u),$
we obtain 
$$
2\ol\pa_+\lan \pa_- u, \,\, \ol\pa_- S\ran_{\tt}=(\ol\pa_+\pa_- u)(\ol\pa_- S)
-(\pa_- u)(\ol\pa_+ \ol\pa_- S).
$$
Since $\ol\pa_+ \ol\pa_- S=0,$ we obtain (\ref{20lpalanpa-uolpa-S}). By the same method, we obtain the result.
\end{proof}
\bgn{proposition}\label{olpadotJpsiolpaS}$\ol\pa\dot{\J_\psi}(\ol\pa S)$ is given by 
\bgn{equation}
\ol\pa\dot{\J_\psi}(\ol\pa S)=8\ol\pa_+\ol\pa_-\lan \pa u, \,\, \J_\psi(\ol\pa S)\ran_{\tt}
\end{equation}
\end{proposition}
\bgn{proof}
Applying Lemma \ref{2olpalanpa-uolpa} to (\ref{dotJpsiolpaS=}), 
we obtain
\bgn{align}
\dot{\J_\psi}(\ol\pa S)=&
8\sqrt{-1}\,\,\ol\pa_-\lan \pa_+u, \,\, \ol\pa_+ S\ran_{\tt}
+8\sqrt{-1}\,\,\ol\pa_+\lan \pa_- u, \,\, \ol\pa_- S\ran_{\tt}
\end{align}
Thus we have 
\bgn{align*}
\ol\pa\dot{\J_\psi}(\ol\pa S)=&8\sqrt{-1}\,\,\ol\pa_+\ol\pa_-\lan \pa_+u, \,\, \ol\pa_+ S\ran_{\tt}
-8\sqrt{-1}\,\,\ol\pa_+\ol\pa_-\lan \pa_- u, \,\, \ol\pa_- S\ran_{\tt}\\
=&8\sqrt{-1}\,\,\ol\pa_+\ol\pa_-\(\lan \pa_+u, \,\, \ol\pa_+ S\ran_{\tt}-\lan \pa_- u, \,\, \ol\pa_- S\ran_{\tt}\)\\
=&8\ol\pa_+\ol\pa_-\(\lan \pa_+u, \,\, \J_\psi(\ol\pa_+ S)\ran_{\tt}
+\lan \pa_- u, \,\, \J_\psi(\ol\pa_- S)\ran_{\tt}\)\\
=&8\ol\pa_+\ol\pa_-\(\lan \pa_+u, \,\, \J_\psi(\ol\pa S)\ran_{\tt}
+\lan \pa_- u, \,\, \J_\psi(\ol\pa S)\ran_{\tt}\)\\
=&8\ol\pa_+\ol\pa_-\lan \pa u, \,\, \J_\psi(\ol\pa S)\ran_{\tt}
\end{align*}
\end{proof}
%%%%%%%%%%%%%%%%%%%%%%%%%%%%%%%%%%%%%%%%%%%%%%%%%

\bgn{proposition}
Let $S_s$ be a family of real functions which smoothly depends a parameter $s $ and satisfies  $S_0=S(\J)$ and 
\bgn{align}\label{dot S=}
\dot S =-4\sqrt{-1}\lan\J_\psi \ol\pa S, \,\, \pa u\ran_{\tt}.
\end{align}
Then the family $\{S_s\}$ satisfies the following 
\bgn{equation}\label{olpaJpsiSs}
\frac{d}{ds}(\ol\pa\J_{\psi_{-s}}\ol\pa) S_s\Big|_{s=0}=0,\end{equation}
where $\J_{\psi_{-s}}$ is the generalized complex structure induced from $\psi_{-s}.$
\end{proposition}

\bgn{proof}
(\ref{olpaJpsiSs}) is written as
\bgn{align*}
0=&\frac{d}{ds}\ol\pa\J_{\psi_{-s}}\ol\pa S_s\Big|_{s=0}\\
=&\ol\pa\dot{\J_{\psi}}(\ol\pa S)+\ol\pa \J_\psi (\ol\pa \dot S)\\
\end{align*}
From $\ol\pa\J_\psi\ol\pa\dot S=-2\sqrt{-1}\,\ol\pa_+\ol\pa_- \dot S$ and Proposition \ref{olpadotJpsiolpaS}, we obtain
\bgn{align}
\ol\pa\dot \J_\psi(\ol\pa S)+\ol\pa\J_\psi(\ol\pa\dot S)=8\ol\pa_+\ol\pa_-\lan \pa u, \,\, \J_\psi(\ol\pa S)\ran_{\tt}-2\sqrt{-1}\ol\pa_+\ol\pa_- \dot S
\end{align}
Since $\dot S$ satisfies 
$$
\dot S=-4\sqrt{-1}\lan \pa u, \,\, \J_\psi(\ol\pa S)\ran_{\tt}
$$
(\ref{olpaJpsiSs}) holds.
\end{proof}
Let $L(\J, \J_{\psi_{-s}})$ be the fourth order differential operator 
$(\ol\pa\J_{\psi_{-s}}\ol\pa)^{*}(\ol\pa\J_{\psi_{-s}}\ol\pa)$, where 
$(\ol\pa\J_{\psi_{-s}}\ol\pa)^{*}$ denotes the adjoint of $(\ol\pa\J_{\psi_{-s}}\ol\pa).$  

\bgn{proposition}\label{fracddsLS=LolL-Lu}
Let $F_s$ be deformations given by $\X_{\sqrt{-1}u_{\Re}}:={\sqrt{-1}}\J_\psi \ol\pa u_{\Re}-\sqrt{-1}\J_\psi\pa u.$
We denote by $L_s$ the operator $L( F_{s\#}\J, \J_\psi).$ We assume that $\ol\pa\J_\psi\ol\pa S=0$. 
Then we have
$$
\(\frac{d}{ds}L_s\)S\Big|_{s=0}= 2L(\ol L- L)u
$$
\end{proposition}

\bgn{proof}

%%%%%%%%%%%%%%%%%%%%%%%%%
Since $\ol\pa\J_\psi\ol\pa S=0$, we have  $L(\J, \J_\psi)S=0$.
We can take a smooth family of functions $\{S_s\}$ which satisfies $S_0=S$ and (\ref{olpaJpsiSs}).
(It is not necessary that $S_s$ arises as scalar curvature.)
Then (\ref{olpaJpsiSs}) which  is equivalent to
\bgn{equation}\label{LJJpsi-s SsBog s=0}
\frac{d}{ds}\(L(\J, \J_{\psi_{-s}}) S_s\)\Big|_{s=0}=0
\end{equation}
%%%%%%%%%%%%%%%%%%%%%%%%%%%
From (\ref{FL(J, Jpsi)}), we have
$$
L(F_{s\#} \J, \J_\psi)f_s^*S_s =
L(F_{s\#} \J, F_{s\#}(F_{-s\#}\psi))f_s^*S_s=
f_s^*\(L(\J,\,\, F_{-s\#}\psi) S_s\).
$$
Then it follows that 
$$
\frac{d}{ds}L(F_{s\#} \J, \J_\psi)f_s^*S_s =\frac{d}{ds}f_s^*\(L(\J, F_{-s\#}\J_\psi) S_s\).
$$
Thus (\ref{LJJpsi-s SsBog s=0}) is equivalent to
\bgn{equation}\label{fracddsLFsJpsi}
\frac{d}{ds}L(F_{s\#} \J, \J_\psi)f_s^*S_s \Big|_{s=0}=0
\end{equation}
We denote by $L_s$ the operator $L( F_{s\#}\J, \J_\psi).$
Then (\ref{fracddsLFsJpsi})
 is equivalent to 
 \bgn{align}\label{fracdds LS}
\(\frac{d}{ds}L_s\)S\Big|_{s=0}+L_s(-(\J e_u)S+\dot S)\Big|_{s=0}=0,
\end{align}
where $-(\J e_u) S=\frac{d}{ds}f_s^* S\Big|_{s=0}=-2\lan \J e_u, dS\ran_{\tt}$.
Since $-\J_\psi(du)=e_u,$ we have
\bgn{align}\label{frac12LJeu}
-\L_{\J(e_u)}S=&-2\lan \J(e_u), \,\,dS\ran_{\tt}
=2\lan\J\J_\psi d u, \,\, dS \ran_{\tt}\\
=&2\lan \J_\psi \J d u , \,\, dS\ran_{\tt}=-2\lan \J d u , \J_\psi dS\ran_{\tt}\\
=&-2\lan \sqrt{-1}\,\ol\pa u-\sqrt{-1}\pa u, \,\, \J_\psi(\pa S+\ol\pa S)\ran_{\tt}\\
=&-2\lan \sqrt{-1}\,\ol\pa u, \,\, \J_\psi \pa S\ran_{\tt}+2\lan \sqrt{-1}\pa u, \,\, \J_\psi\ol\pa S\ran_{\tt}
\end{align}

Then from (\ref{fracdds LS}), (\ref{dot S=}) and (\ref{frac12LJeu}), we obtain
\bgn{align*}
\(\frac{d}{ds}L_s\)S\Big|_{s=0}=&-L_s(-\L_{\J (e_u)}S+\dot S)\Big|_{s=0}\\
=-&L\(-2\lan \sqrt{-1}\,\,\ol\pa u, \,\, \J_\psi \pa S\ran_{\tt}+2\lan \sqrt{-1}\pa u, \,\, \J_\psi\ol\pa S\ran_{\tt}\)\\
-&L\(-4\sqrt{-1}\lan \pa u,\,\,\J_\psi \ol\pa S\ran_{\tt}\)\\
=&2L\(\lan \sqrt{-1}\,\,\ol\pa u, \,\, \J_\psi \pa S\ran_{\tt}+\lan \sqrt{-1}\pa u, \,\, \J_\psi\ol\pa S\ran_{\tt}\)\\
=&2\sqrt{-1}L\(\lan du, \,\, \J_\psi dS\ran_{\tt}\)\\
=-&2\sqrt{-1}L\lan \J_\psi du,\,\, dS\ran_{\tt}\\
=&-\sqrt{-1}L\{ u, S\}_{\J_\psi}
\end{align*}
Applying Proposition \ref{-iLolL=isJ}, i.e.,
$\frac2{\sqrt{-1}}(L-\ol L)u =\{u, S(\J)\}_{\J_\psi}$, we obtain

$$
\(\frac{d}{ds}L_s\)S\Big|_{s=0}=2L(\ol L-L)S
$$
\end{proof}
\bgn{proposition}\label{LolL=lolL}
Let $(M, \J, \J_\psi)$ be a generalized extremal K\"ahler manifold. 
Then we have
$$L\ol L =\ol L L$$
\end{proposition}
\bgn{proof}
Two pure imaginary functions $\sqrt{-1}u_1, \sqrt{-1}u_2$ gives 
$\X_{\sqrt{-1}u_1}$ and $\X_{\sqrt{-1}u_2},$ respectively. 
Then $\X_{\sqrt{-1}u_1}$ and $\X_{\sqrt{-1}u_2}$ gives rise to 
$2$-parameter deformations $\{\J_{t_1, t_2}\}.$
We shall calculate the Hessian of the functional 
$\E(\J):=\int_M S(\J)^2\vol_M$ under the deformations $\{\J_{t_1, t_2}\}.$
From Proposition \ref{Jtu|t=0=L+olL u},
the differential of $\Phi(\J_{t_1, t_2})$ with respect to $t_1$ is given by
\bgn{align}
\frac12\frac{\pa}{\pa t_1}\Phi(\J_{t_1, t_2})=&\frac12\int_M S(\J_{t_1, t_2})\frac{d}{d t_1}S(\J_{t_1,t_2})\Big|_{t_1=0}\\
=&\int_M S(\J_{t_1, t_2})\((L_{t_1,t_2}+\ol L_{t_1, t_2})u_1\)\vol_M\\
=&\int_M u_1\((L_{t_1,t_2}+\ol L_{t_1, t_2})S(\J_{t_1, t_2})\)\vol_M
\end{align}
From Proposition \ref{-iLolL=isJ}, we also have 
$$
(L_{t_1,t_2}-\ol L_{t_1, t_2})S(\J_{t_1, t_2})=\{S(\J_{t_1, t_2}),\,\,S(\J_{t_1, t_2})\}_{\J_\psi}=0
$$
Thus 
\bgn{equation}\label{fracpapat1PhiJt1t2=2intM}
\frac12\frac{\pa}{\pa t_1}\Phi(\J_{t_1, t_2})=2\int_M u_1 L_{t_1, t_2}S(\J_{t_1, t_2})=2\int_M u_1\ol L_{t_1, t_2}S(\J_{t_1, t_2})
\end{equation}
From Proposition \ref{fracddsLS=LolL-Lu}, i.e., $\dot L S=2L(\ol L- L)u$ and Proposition \ref{Jtu|t=0=L+olL u}, i.e.,
$\frac{d}{dt}S(\J_t^{u})|_{t=0}=2(L+\ol L)u_{\Im}$,
the differential of $\frac{\pa}{\pa t_1}\Phi $ with respect to $t_2$ is given by
\bgn{align}
\frac12\frac{\pa}{\pa t_2}\frac{\pa }{\pa t_1}\Phi(\J_{t_1,t_2})\Big|_{t_1, t_2=0}=&
2\int_M u_1\frac{\pa}{\pa t_2}L_{t_1, t_2}S(\J_{t_1, t_2})\Big|_{t_1, t_2=0}\\
=&2\int_M u_1\dot L_{t_1, t_2}S(\J_{t_1, t_2})\Big|_{t_1, t_2=0}\\
+&2\int_M
u_1L_{t_1, t_2}\dot S(\J_{t_1, t_2})\Big|_{t_1, t_2=0}\\
=&4\int_M u_1L(\ol L-L) u_2+u_1L(L+\ol L) u_2\\
=&8\int_M u_1(L\ol L u_2 )\,\,\vol_M
\end{align}
From (\ref{fracpapat1PhiJt1t2=2intM}),
the similar calculation gives 
\bgn{align}
\frac12\frac{\pa}{\pa t_2}\frac{\pa}{\pa t_1}\Phi(\J_{t_1,t_2})\Big|_{t_1, t_2=0}=&
2\int_M u_1\frac{d}{d t_2}\ol L_{t_1, t_2}S(\J_{t_1, t_2})\Big|_{t_1, t_2=0}\\
=&4\int_M u_1\ol L( L-\ol L) u_2+u_1\ol L(L+\ol L) u_2\\
=&8\int_M u_1(\ol L L u_2) \,\,\vol_M
\end{align}
Hence we obtain $L\ol L u=\ol L L u$.
\end{proof}
Let $(M, \J, \J_\psi)$ be a generalized K\"ahler manifold with generalized scalar curvature $S(\J).$
The generalized metric $G=-\J\circ\J_\psi$ defines the generalized isometry group $I_G(M)$ which is the subgroup of $\wtil\Diff(M)$ preserving $G.$
Then it turns out that $I_G(M)$ is a compact Lie group.
The Lie algebra of $I_G(M)$ is denoted by $\frak i_G(M)$.
If $(M,\J, \J_\psi)$ is a generalized extremal K\"ahler manifold, then 
the scalar curvature $S:=S(\J)$ gives the class $[\J_\psi\circ \ol\pa S]\in H^1(\w^\bullet\ol\L_\J)$ and the adjoint action of $[\J_\psi \ol\pa S]$ on $\frak g_0$ gives the decomposition into $\lam$-eigenspaces 

$$
{\frak g}({\lam}):= \{ a\in {\frak g}_{\red}\,|\, \ad_{[\J_\psi\ol\pa S]}\a=\lam \a\,\,\}$$

\bgn{theorem} Let $(M, \J, \J_\psi)$ be a generalized extremal K\"ahler manifold with generalized 
scalar curvature $S(\J).$
Then the Lie algebra of the reduced automorphisms $\frak g_{\red}$ of $(M,\J,\J_\psi)$ admits the following decomposition as Lie algebra:

$$
{\frak g}_{\red}={\frak g}(0)\oplus \sum_{\lam\neq 0}{\frak g}(\lam),
$$
where ${\frak g}(0)$
is the maximal reductive subalgebra $(\frak i_G(M)\cap {\frak g}_{\red})\otimes\C. $
\end{theorem}
\bgn{proof}
Since the Lie algebra ${\frak g}_{\red}$ of the reduced automorphisms is identified with the space of complex functions which are annihilated by the action of $L,$
$$
\{\, u \in C_0^\infty(M,\C)\, |\, Lu=0\, \},
$$
where $Lu=(\ol\pa\J_\psi\ol\pa)^*(\ol\pa\J_\psi\ol\pa)u.$
From Proposition \ref{LolL=lolL}, we have $L\ol L=\ol L L.$
Thus the action of $\ol L$ preserves the kernel space $\{\, u \in C_0^\infty(M,\C)\, |\, Lu=0\, \}$ of $L$
and then we have the eigenspace decomposition of the action of $\ol L$ under the identification, 
$${\frak g}_{\red} =\oplus_{\lam}V_{-\sqrt{-1}\lam}
$$
Then from Proposition \ref{-iLolL=isJ}, i.e.,
$\sqrt{-1}(L-\ol L)u =\{u, S(\J)\}_{Poi}$, it follows that $-\sqrt{-1}\lam$-eigenfunction $u$ satisfies the following: 
\bgn{align}
-\sqrt{-1}\lam u =&\ol L u =(\ol L-L)u \\
=&-\sqrt{-1}\{ S(\J), \,\, u\}_{Poi} 
\end{align} 
Thus we have 
$$
[\J_\psi\ol\pa S(\J), \,\, \J_\psi\ol\pa u]_{\cou}=\{ S(\J), \,\, u\}_{Poi} =\lam u.
$$
Hence we have 
${\frak g}(\lam) =V_{-\sqrt{-1}\lam}.$
Since $V_0=\ker L\cap \ker\ol L,$ we have 
$V_0=(\frak i_G(M)\cap {\frak g}_{\red})\otimes\C.$
Thus we have the result.
\end{proof}

%%%%%%%%%%%%%%%%%%%%%%%%%%%%%%%%

\medskip
\noindent
E-mail address: goto@math.sci.osaka-u.ac.jp\\
\noindent
Department of Mathematics, Graduate School of Science,\\
\noindent Osaka University Toyonaka, Osaka 560-0043, JAPAN

\end{document}